\documentclass{article}
\usepackage[utf8]{inputenc}
\usepackage[margin=1in]{geometry}
\usepackage{ dsfont }
\usepackage{microtype}
\usepackage{graphicx}
\usepackage{subfigure}
\usepackage{booktabs}
\usepackage{hyperref}
\usepackage{amsmath, amsfonts, amssymb, amsthm}
\usepackage{xcolor}
\usepackage{algorithm}
\usepackage{algorithmicx}
\usepackage{algpseudocode}
\usepackage{booktabs}
\usepackage{graphicx}
\usepackage{bbm}
\usepackage{bm}
\usepackage{authblk}
\usepackage{thmtools}
\usepackage[round]{natbib}
\newtheorem{theorem}{Theorem}
\newtheorem{lemma}[theorem]{Lemma}
\newtheorem{proposition}[theorem]{Proposition}
\newtheorem{remark}[theorem]{Remark}

\newcommand{\blockcomment}[1]{}

\newcommand{\E}{\mathbb{E}}
\newcommand{\V}{\mathbb{V}}
\newcommand{\R}{\mathbb{R}}
\renewcommand{\P}{\mathbb{P}}

\newcommand{\Z}{\mathcal{Z}}

\renewcommand{\a}{\alpha}
\renewcommand{\b}{\beta}
\renewcommand{\d}{\delta}

\newcommand{\g}{\gamma}

\newcommand{\s}{\sigma}
\DeclareMathOperator*{\argmax}{argmax}

\def\Z{\mathbb{Z}}

\def\mE{\mathcal{E}}
\def\mA{\mathcal{A}}

\def\mN{\mathcal{N}}
\def\mH{\mathcal{H}}
\def\mB{\mathcal{B}}
\def\P{\mathbb{P}}
\def\l{\left}
\def\r{\right}
\def\la{\left\langle}
\def\ra{\right\rangle}
\def\h{\hat}
\def\t{\tilde}
\def\pt{\partial}

\def\dt{\d_t}
\def\bs{{\bf s}}
\def\bq{{\bf q}}
\def\bi{{\bf i}}
\def\bj{{\bf j}}
\def\ba{{\bf a}}
\def\be{{\bf e}}
\def\bba{{\bm \a}}
\def\bpi{{\bm \pi}}
\def\bnu{{\bm \nu}}
\def\hs{\h{\bs}}
\def\hq{\h{\bq}}
\def\hls{\h{s}}
\def\hlq{\h{q}}
\def\hc{\h{c}}
\def\hb{\h{\b}}
\def\ha{\h{\a}}
\def\S{\Sigma}

\def\hpi{\h{\bpi}}
\def\tv{\t{v}}
\def\tcb{\textcolor{black}}

\title{Continuous-in-time Limit for Bayesian Bandits}
\author{Yuhua Zhu\thanks{Department of Mathematics and Halıcıoğlu Data Science Institute, University of California-San Diego. (yuz244@ucsd.edu)} \quad
Zach Izzo\thanks{Department of Mathematics, Stanford University. (zizzo@stanford.edu)}
\quad
Lexing Ying\thanks{Department of Mathematics and Institute for Computational and Mathematical Engineering, Stanford University. (lexing@stanford.edu)}
}
\date{}

\begin{document}

\maketitle

\begin{abstract}
This paper revisits the bandit problem in the Bayesian setting. The Bayesian approach formulates the bandit problem as an optimization problem, and the goal is to find the optimal policy which minimizes the Bayesian regret. 
One of the main challenges facing the Bayesian approach is that computation of the optimal policy is often intractable, especially when the length of the problem horizon or the number of arms is large. In this paper, we first show that under a suitable rescaling, the Bayesian bandit problem converges toward a continuous Hamilton-Jacobi-Bellman (HJB) equation. The optimal policy for the limiting HJB equation can be explicitly obtained for several common bandit problems, and we give numerical methods to solve the HJB equation when an explicit solution is not available. Based on these results, we propose an approximate Bayes-optimal policy for solving Bayesian bandit problems with large horizons. Our method has the added benefit that its computational cost does not increase as the horizon increases. \end{abstract}

\section{Introduction}

Bandit problems were first introduced by \cite{thompson1933likelihood} with later pioneering work due to \cite{robbins1952some} and \cite{wald2004sequential}. In more recent years, bandit algorithms have become widely adopted for automated decision-making tasks such as dynamic pricing \citep{ferreira2018online}, mobile health \citep{tewari2017ads}, Alpha Go \citep{silver2016mastering},  etc.

The bandit problem can be considered from one of two perspectives: Bayesian or frequentist. The Bayesian approach dominated bandit research from 1960-1980 \citep{bradt1956sequential,gittins1979bandit}. The objective is to minimize an average cumulative regret with respect to the Bayesian prior measure of the problem environment. It formulates the bandit problem as an optimization problem, and the goal is to find the optimal policy which minimizes the Bayesian regret. In the frequentist setting \citep{lai1985asymptotically}, the cumulative regret is viewed as an unknown deterministic quantity, and the goal is to design policies to achieve the best environment-dependent performance. 

The main difficulty with Bayesian bandits is that computation of the optimal policy is often intractable, especially when the number of arms or the horizon is large. Gittin's index \citep{gittins1979bandit} reduced the computational cost for the discounted infinite horizon setting but does not apply to undiscounted cases \citep{berry1985bandit}. Although computing the Bayes-optimal policy is challenging, there is a significant payoff: the performance of the policy is not only optimal in the Bayesian setting (by definition) but also has favorable frequentist regret guarantees \citep{lattimore2016regret}. In addition, Bayesian bandits have been widely employed in economics \citep{bergemann2006bandit}, including in contract theory \citep{gur2022sequential}, dynamic pricing \citep{leloup2001dynamic}, portfolio management \citep{el2005optimal}, etc.

In this paper, we revisit the Bayesian perspective for the multi-armed bandit problem and analyze it using tools from PDEs. A continuous-in-time limiting HJB equation is derived as the horizon $n$ goes to infinity for a range of bandit problems. Based on the limiting equation, a regularized Bayes-optimal policy is proposed, where regularization is employed to increase exploration and stability. Numerical schemes can be used to approximate the optimal policy, leading to improved computational efficiency when the horizon is large. In addition, the exact optimal policy for the limiting HJB equation can be obtained for certain types of bandit problems, including the classical Bernoulli and Gaussian arm reward cases, resulting in an efficient algorithm to approximate the optimal policy even if the number of arms is large.
In summary, our contributions are as follows.
\begin{itemize}
    \item We derive a continuous-in-time limit, an HJB equation,  for the Bayesian multi-armed bandit problem.
    \item We propose a regularized version of the Bayesian multi-armed bandit problem, which encourages exploration and smooths the optimal policy.
    \item Based on the limiting PDE, we give an efficient algorithm for approximating the optimal policy.
\end{itemize}

Recently, the use of differential equations to analyze machine learning algorithms has received growing interest, especially in optimization algorithms  \cite{su2014differential, li2017stochastic}, sampling algorithms \cite{welling2011bayesian, liu2017stein}, neural networks \cite{mei2018mean,chen2018neural, weinan2018mean,chen2020dynamical}. However, fewer connections have been built between multi-armed bandits and differential equations until the recent two years. For bandit problems, \cite{fan2021diffusion,wager2021diffusion, kobzar2022pde} model several policies  in the frequentist setting via a continuous SDE or PDE, and this continuous analysis is used to provide insights into the properties of the algorithms studied. In the Bayesian setting, \cite{araman2022diffusion, che2018recommender} give differential equation approximations, but their work can only be applied to settings with two possible environments. In this paper, we consider Bayesian bandits with general environments. There are also related works using differential equations in online learning settings, including contextual bandits \citep{kapralov2011prediction}, drifting games \citep{wang2022new}, etc.

\section{Bayesian Bandits}
Throughout the paper, we focus on $K$-armed stochastic bandits played over $n$ rounds, where $n \in \Z_+$ is a positive integer called the horizon. At each round $i$, the learner chooses an arm (also called an action) $A^i$ from the action space $\mA = \{a_k\}_{k=1}^K$ according to a policy $\pi^i$, and the environment reveals a reward $X^i\in\R$. The underlying environment $\nu$ belongs to an environment class $\mE$ and defines the arm reward distributions. More precisely, given an environment $\nu = (P^\nu_a, a\in\mA)\in\mE$, the reward $X^i$ follows the distribution $P_{A^i}^\nu$.

The policy $\bpi^i$ at round $i$ is a function which maps the history $H^{i} = (A^1,X^1, \cdots, A^{i-1},X^{i-1})$ to a probability distribution over the action space $\mA$. More precisely, we denote the set of possible histories at the beginning of round $i$ by $\mH^i = (\mA\times\R)^{i-1}$, and $\mH^1 = \emptyset$. We denote by $\Delta(\mA)$ the set of probability measures on $\mA$ so that the policy $\bpi^i$ is a mapping from $\mH^i$ to $\Delta(\mA)$. We denote by $\Pi$ the set of policies $\bpi = \{\bpi^i\}_{i=1}^n$, which is measurable with respect to the filtration associated with the process $\{H^i\}_{i=1}^n$. We call $\Pi$ the competitor class. 

Define the expectation of arm $a$ in environment $\nu$ as $\mu_a(\nu)$, where
\begin{equation}\label{eq:defmua}
    \mu_a(\nu) = \int_\R x P^\nu_a(x)dx.
\end{equation}
The expected cumulative reward $c_n(\bpi,\nu)$ measures the performance of policy $\bpi$ in environment $\nu$,
\[
c_n(\bpi,\nu) = \E\l[\sum_{i=1}^n\mu_{A^i}(\nu)\r],
\]
where the expectation is taken over the probability measure induced by the interaction of the policy and the environment. 
The goal in the $K$-armed bandit setting is to design a policy $\bpi^*$ that leads to the largest expected cumulative reward among all policies in the competitor class $\Pi$. The main difficulty arises because the environment is unknown, and the policy can only depend on the history sequences $\{H^i\}_{i=0}^{n-1}$. 

One way to measure the performance of a policy $\bpi$ is to find functions $C:\mE\to[0,\infty), f:\mE\to[0,\infty)$ that upper bound the regret:
\[
R_n(\bpi,\nu) = \tcb{n}\mu^*(\nu) - c_n(\bpi,\nu) \leq C(\nu)f(n),
\]
where $\mu^*(\nu) = \max_a\mu_a(\nu)$ is the expected reward of the optimal arm. This is the frequentist regret, which is environment-dependent \citep{lattimore2020bandit}.

\tcb{Another way to measure the performance of a policy is via the averaged cumulative reward with respect to a probability measure $\rho(\nu)$ on the environment $\mE$}, 
\[
c_n(\bpi,\rho) = \E\l[ \sum_{i=1}^n\int_\mE\mu_{A^i}(\nu)\rho(\nu)d\nu \r].
\]
In the Bayesian setting, the environment is viewed as a random variable. According to Bayes' rule, the probability $\rho(\nu)$ will be updated conditional on the history sequence. Given a horizon $n$, we assume that at each round $i$, the environment $\nu^i$ is sampled from a prior measure $\rho^i(\nu)$ over the environment class $\mE$. After pulling arm $A^i$ and obtaining the reward $X^i\sim P_{\nu}^{A^i}$, we update $\rho^{i+1}(\nu)$ to be the posterior distribution of the environment. Given an initial prior measure $\rho^1(\nu)$, the goal is to find the optimal policy that maximizes the averaged cumulative reward,
\[
\begin{aligned}
\max_{\bpi\in\Pi} c_n(\bpi) = \E\l[\sum_{i=1}^n\int_{\mE}\mu_{A^i}(\nu) \rho^i(\nu)d\nu\r].
\end{aligned}
\]
The bandit problem with the above objective function is called a {\em Bayesian bandit} (Chapter 35 in \cite{lattimore2020bandit}).



\color{black}
\section{Continuous Limits of Bayesian Bandits}
\subsection{An illustrative example}\label{sec:example}
Consider the ``one-armed'' bandit problem in which the reward of the first arm follows a Bernoulli$(\nu)$ distribution (with $\nu$ unknown) and the second arm gives a deterministic reward $\mu_2 = \frac12$. We assume an initial prior distribution of $\nu\sim\text{Beta}(\a,\b)$. Then the posterior measure of $\nu$ at round $i$ depends on two quantities. The first quantity is $q^i$, which is the number of pulls of the unknown arm before round $i$. The second quantity is $s^i$, which is the cumulative reward of the unknown arm before round $i$. The posterior distribution of $\nu$ is $\text{Beta}(\a+s^i,\b+q^i-s^i)$. 

Let $w^i(s,q)$ be the optimal cumulati\color{black}ve reward starting from round $i$ with $s^i = s$, $q^i = q$.
\color{black}
Similarly, let $w^i_k(s,q)$ be the optimal cumulative reward in this same setting, assuming that the $k$-th arm is pulled at round $i$. Formally, we have
$$ w^i(s, q) = \max_{\bpi\in\Pi} \E\l[\sum_{j=i}^n\int_{\mE}\mu_{A^j}(\nu) \rho^j(\nu)d\nu \: \Bigg| \: s^i = s, q^i = q \r], $$
$$ w^i_k(s, q) = \max_{\bpi\in\Pi} \E\l[\sum_{j=i}^n\int_{\mE}\mu_{A^j}(\nu) \rho^j(\nu)d\nu \: \Bigg| \: s^i = s, q^i = q, A^i = k \r]. $$
As before, the dependence of the expectation on the policy $\bpi$ is through the actions $A^j$.
\color{black}
If the second arm is chosen at round $i$, then 
\[
w_2^i(s,q) = \mu_2 + w^{i+1}(s,q).
\]
If the first arm is chosen at round $i$, then 
\[
w_1^i(s,q) = p(s,q)+p(s,q)w^{i+1}(s+1,q+1)+(1-p(s,q))w^{i+1}(s,q+1),
\]
where \[
p(s,q) = \frac{\a+s}{\a+\b+q}.
\]
The first term $p(s,q)$ represents the expectation of the reward if the first arm  is pulled at round $i$. The second and third term hold because after pulling the first arm, one has
\[
\begin{aligned}
    &q^{i+1} = q^i+1;\\
    &\P(s^{i+1} = s^i+1|s^i) = p(s^i,q^i),
    &\P(s^{i+1} = s^i|s^i) = 1-p(s^i,q^i).
\end{aligned}
\]
Since $w^i$ is the optimal reward, one has
\begin{equation}\label{eq:bdt-algo}
    w^i(s,q) = \max\{w^i_1(s,q), w^i_2(s,q)\}.
\end{equation}
Note that for horizon $n$, $w^{n+1}(s,q) = 0$ for all $s, q$ by the definition of $w$. Therefore, one can compute $w^i(s,q)$ for all $s=\{0,\cdots,i-1\}, q = \{0,\cdots,i-1\}$ via backwards induction.

\begin{figure}[h!]
    \centering
    \includegraphics[width=0.6\linewidth]{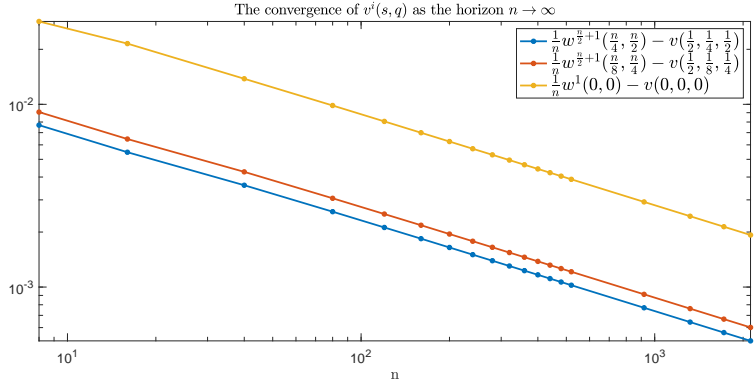}
    \caption{The above plot shows the decrease in the error $|\frac{1}{n}w^i(s,q) - v(t,\hls,\hlq)|$ as the horizon $n\to\infty$ for $(i,s,q) = \{(\frac{n}{2}+1, \frac{n}{4},\frac{n}{2}),(\frac{n}{2}+1, \frac{n}{8},\frac{n}{4}),(1, 0,0)\}$ and the corresponding $(t,s,q) = \{(\frac12,\frac14,\frac12)(\frac12,\frac18,\frac14),(0,0,0)\}$. We set the initial hyperparameters $(\a,\b) = (\frac n2,\frac n2)$,.}
    \label{fig:convergence}
\end{figure}

To derive a continuous limit in this setting, we rescale the reward and the number of arm pulls by $1/n$:
\[
\hls = \frac1ns,\quad \hlq = \frac1n.
\]
One can define the rescaled optimal reward function
\[
v^i(\hls,\hlq) = \frac1nw^i(s,q),
\] which then satisfies the equation
\begin{equation*}
\begin{aligned}
    v^i(\hls,\hlq) = \max&\l\{\frac1n\mu_2 + v^{i+1}(\hls,\hlq), \r.\\
    &\l.\frac1n\t{p}(\hls,\hlq)+\t{p}(\hls,\hlq)v^{i+1}(\hls+\frac1n,\hlq+\frac1n)+(1-\t{p}(\hls,\hlq))v^{i+1}(\hls,\hlq+\frac1n)\r\},    
\end{aligned}
\end{equation*}
where
\begin{equation*}
    \t{p}(s,q) = \frac{n^{-1}\a + s}{n^{-1}(\a+\b)+q}.
\end{equation*}
By setting $\d_t = \d_s = \d_q = n^{-1}$, the above equation can be equivalently written as
\begin{equation*}
\begin{aligned}
    \frac{v^{i+1}(\hls,\hlq) - v^i(\hls,\hlq)}{\d_t} + \max&\l\{\mu_2, \ \t{p}(\hls,\hlq)+\t{p}(\hls,\hlq)\frac{v^{i+1}(\hls+\d_s,\hlq+\d_q)-v^{i+1}(\hls,\hlq+\d_q)}{\d_s} \r.\\
    &\l.\quad\quad\quad+ \frac{v^{i+1}(\hls,\hlq+\d_q) - v^{i+1}(\hls,\hlq)}{\d_q}\r\} = 0.    
\end{aligned}
\end{equation*}
From the above equation, one sees that the rescaled value function $v^i(\hls,\hlq)$ is the numerical discretization of the following PDE:
\begin{equation}\label{eq:bdt-hjb1}
    \pt_tv(t,\hls,\hlq) + \max\{\mu_2, \h{\mu}(\hls,\hlq) + \h{\mu}(\hls,\hlq)\pt_{\hls}v(t,\hls,\hlq) + \pt_{\hlq}v(t,\hls,\hlq)\} = 0,\quad v(1,\hls,\hlq) = 0,
\end{equation}
where 
\begin{equation*}
    \h{\mu}(s,q) = \lim_{n\to\infty}\t{p}(s,q).
\end{equation*}
Furthermore, by moving the constant $\mu_2$ outside of the maximum operator and introducing a control parameter $\h{\pi}\in[0,1]$, one arrives at a Hamilton-Jacobi-Bellman equation for $v$:
\begin{equation}\label{eq:bdt-hjb}
    \pt_tv + \max_{\h{\pi}(t,\hls,\hlq)\in[0,1]}\l(\h{\mu} + \h{\mu}\pt_{\hls}v + \pt_{\hlq}v-\mu_2\r)\pi + \mu_2 = 0, \quad \quad v(1,\hls,\hlq) = 0.
\end{equation}
\eqref{eq:bdt-hjb1} and \eqref{eq:bdt-hjb} are equivalent because when $\h{\mu} + \h{\mu}\pt_{\hls}v + \pt_{\hlq}v > (<) \mu_2$, then $\h{\pi} = 1 (=0)$, respectively. 
In other words, as the horizon $n\to\infty$, the rescaled value function $v^i(\hls,\hlq)$ satisfies the above HJB equation. We plot the convergence of $\frac{1}{n}w^i(\hls,\hlq)$ as $n$ increases in Figure~\ref{fig:convergence}.

Classical results from optimal control \tcb{(see, e.g., Chapter 10.3.3 of \cite{evans2010partial})} imply that $v(t,\hls,\hlq)$ in \eqref{eq:bdt-hjb} solves the  control problem
\begin{equation}\label{eq:bdt-control}
\begin{aligned}
    \max_{\h{\pi}(\tau)\in[0,1]}\quad &\int_t^1 \h{\mu}(\hls(\tau),\hlq(\tau))\pi(\tau) + (1-\pi(\tau))\mu_2 d\tau\\
    \text{s.t.} \quad &d\hlq(\tau) = \pi(\tau) d\tau,\\
    &d\hls(\tau) = \h{\mu}(\hls(\tau),\hlq(\tau))\pi(\tau) d\tau,\\
    &\hlq(t) = \hls, \quad \hls(t) = \hlq.
\end{aligned}
\end{equation}

The preceding example illustrates that the Bayesian bandit algorithm \eqref{eq:bdt-algo} can be viewed as the discretization of an HJB equation \eqref{eq:bdt-hjb}, which solves the control problem \eqref{eq:bdt-control}. In other words, as the horizon $n\to\infty$, the Bayesian bandit problem will converge to a continuous control problem that can be solved via the HJB equation. In the next section, we extend this formulation to a more general setting.

\subsection{\tcb{Formal derivation} from Bayesian bandits to the HJB equation}
We return to the original $K$-armed bandit setting with horizon $n$ and environments parameterized by $\bnu\in\mE\subset \R^{d}$. The reward of the $k$-th arm $a_k$ follows the distribution $P_k^\bnu$.  
We assume a prior measure $\rho(\bnu)$ on the environment $\bnu$ at the beginning of the first round. We also assume that the updated measure $\rho^i(\bnu)$ at the beginning of round $i$ only depends on $(\bs^i,\bq^i)$. Here $\bs^i = (s^i_k)_{k=1}^K\in\R^K$ is a $K$-dimensional vector representing the cumulative reward of each arm up to round $i-1$, and $\bq^i = (q^i_k)_{k=1}^K\in\{0,\cdots,i-1\}^K$ is a $K$-dimensional vector representing the number of pulls of each arm up to round $i-1$.
\color{black}
We note that the assumption that the state space can be reduced to $(\bs^i,\bq^i)$ covers many, but not all, bandit algorithms. For instance, many stochastic bandit algorithms, where the arm rewards are drawn from a stationary probability distribution, can be represented using this state space, but the Exp3 algorithm \citep{auer2002exp3} cannot. For an extensive discussion, refer to Section~2.1 of \cite{wager2021diffusion}.
\color{black}

Under the above assumptions, the posterior distribution of the environment at round $i$ can be written as a function of $\bs^i,\bq^i$, i.e., $\rho(\bnu|\bs^i,\bq^i)$. If arm $k$ is pulled at round $i$, then the reward $X^i$ follows the distribution $P_k^\bnu$, and 
\begin{equation}\label{eq:sip1}
    s_k^{i+1} = s_k^i+X^i, q^{i+1}_k = q^i_k+1. 
\end{equation}

Let $w^i(\bs,\bq)$ be the optimal expected cumulative reward starting from round $i$ with $(\bs^i,\bq^i) = (\bs,\bq)$. Then it satisfies the following equation:
\begin{equation}\label{eq:wi}
    w^i(\bs,\bq) = \max_k\l\{ \int_\mE \mu_k(\bnu) \rho(\bnu|\bs,\bq) d\bnu + \int_\mE\int_\R w^{i+1}(\bs+x\be_k, \bq+\be_k) P^\bnu_k(x)\rho(\bnu|\bs,\bq)\,dx\, d\bnu \r\},
\end{equation}
where $\mu_k(\bnu)$ is the expected reward of the $k$-th arm defined in \eqref{eq:defmua}, and $\be_k$ is a $K$-dimensional vector with the $k$-th element being $1$ and all other elements being $0$. The first term in the max operator represents the expectation of the rewards if the $k$-th arm is pulled. The second term is due to the fact that $s^{i+1}_k$ and $q^{i+1}_k$ will follow \eqref{eq:sip1} if the $k$-th arm is pulled at round $i$.

Next, we rescale the parameters to derive the continuous-in-time limit. Let 
\begin{equation}\label{eqn: rescale}
    t = \frac{i-1}n, \quad \hq = \frac1n\bq^i,\quad \hs = \frac{1}{f(n)}\bs^i, \quad v\left(\frac{i-1}{n},\frac{\bs}{f(n)}, \frac{\bq}{n}\right) = \frac{1}{f(n)}w^i(\bs,\bq).
\end{equation}
We shrink the $n$ rounds to the time interval $[0,1]$ so that as $n\to\infty$, the discrete round $i\leq n$ will correspond to a continuous time $t=(i-1)/n\in[0,1]$. We also rescale the number of pulls to $[0,1]$, so that it is on the same scale as the rescaled rounds. The cumulative reward $\bs^i$ and $w^i$ are rescaled by $\frac{1}{f(n)}$, where $f(n)$ will be determined later. Accordingly, the rescaled expected reward $v(t,\hs,\hq)$ becomes a function of the continuous time $t$ and rescaled history $(\hs,\hq)$. We refer to the function $f(n)$ as the \emph{scaling factor}. For different scaling factors, the limiting rescaled cumulative reward $v(t,\hs,\hq)$ will follow different dynamics. 

Define the moments of the $k$-th arm w.r.t. the probability measure $P^\bnu_k(x)$ and Bayesian measure $\rho(\bnu|\bs, \bq)$ by
\begin{equation}\label{eq:defofbar}
\begin{aligned}
        &\bar{\mu}_k(\bs, \bq) = \int_\mE\int_\R xP^\bnu_k(x) \rho(\bnu|\bs, \bq)\,dx\,d\bnu, \quad \bar{\sigma}_k^2(\bs, \bq) = \int_\mE\int_\R x^2P^\bnu_k(x) \rho(\bnu|\bs, \bq)\,dx\,d\bnu, \\
        &\bar{E}_k^p(\bs, \bq) = \int_\mE\int_\R x^pP^\bnu_k(x) \rho(\bnu|\bs, \bq)\,dx\,d\bnu.
\end{aligned}
\end{equation}
\color{black}
We assume that these moments exist and are finite.
\color{black}
Inserting the Taylor expansion 
\[
w^{i+1}(\bs+x\be_k, \bq+\be_k) = \sum_{p=0}^\infty\frac{1}{p!}\pt^p_{s_k}w^{i+1}(\bs, \bq+\be_k)x^p
\]
into \eqref{eq:wi} yields
\[
\begin{aligned}
    w^i(\bs,\bq) = \max_k&\l\{ \bar{\mu}_k(\bs,\bq) + w^{i+1}(\bs,\bq+\be_k) + \bar{\mu}_k(\bs,\bq)\pt_{s_k}w^{i+1}(\bs,\bq+\be_k) + \frac12\bar{\sigma}^2_k(\bs,\bq) \pt_{s_k}^2w^{i+1}(\bs,\bq+\be_k)
\r.\\
&\quad \l.+ \sum_{p=3}^\infty\frac{1}{p!}\bar{E}^p_k(\bs,\bq)\pt_{s_k}^pw^{i+1}(\bs,\bq+\be_k) \r\}.
\end{aligned}
\]
Therefore, the rescaled reward $v\left(\frac{i-1}{n},\frac{\bs}{f(n)}, \frac{\bq}{n}\right) = \frac{1}{f(n)}w^i(\bs,\bq)$ satisfies
\[
\begin{aligned}
v\l(\frac{i-1}{n},\frac{\bs}{f(n)},\frac{\bq}{n}\r) = \max_k&\l\{ \frac{1}{f(n)}\bar{\mu}_k(\bs,\bq) + v\l(\frac{i}{n},\frac{\bs}{f(n)},\frac{\bq+\be_k}{n}\r) \r.\\
&\quad+ \frac{1}{f(n)}\bar{\mu}_k(\bs,\bq)\pt_{\h{s}_k}v\l(\frac{i}{n},\frac{\bs}{f(n)},\frac{\bq+\be_k}{n}\r) \\
&\quad + \frac12\frac{1}{f^2(n)}\bar{\sigma}^2_k(\bs,\bq) \pt^2_{\h{s}_k}v\l(\frac{i}{n},\frac{\bs}{f(n)},\frac{\bq+\be_k}{n}\r)\\
&\quad\l.+ \sum_{p=3}^\infty\frac{1}{p!}\frac1{f^p(n)}\bar{E}^p_k(\bs,\bq)\pt^p_{\h{s}_k}v\l(\frac{i}{n},\frac{\bs}{f(n)},\frac{\bq+\be_k}{n}\r)\r\}.
\end{aligned}
\]
After reorganizing the terms and setting $\d_t = \d_q = \frac1n$, one has
\[
\begin{aligned}
\frac{v(t+\d_t,\hs,\hq)-v(t,\hs,\hq)}{\d_t} +\max_k&\l\{ \frac{1}{\d_tf(n)}\bar{\mu}_k(f(n)\hs,n\hq) + \frac{v(t+\d_t,\hs,\hq+\d_q\be_k)-v(t+\d_t,\hs,\hq)}{\d_q} \r.\\
&\l.  + \frac{1}{\d_tf(n)}\bar{\mu}_k(f(n)\hs,n\hq)\pt_{\h{s}_k}v(t+\d_t,\hs,\hq+\d_q\be_k) \r.\\
&\l. + \frac12\frac{1}{\d_tf^2(n)}\bar{\sigma}^2_k(f(n)\hs,n\hq) \pt^2_{\h{s}_k}v(t+\d_t,\hs,\hq+\d_q\be_k)\r.\\
&\l. + \sum_{p=3}^\infty\frac{1}{p!}\frac1{\d_tf^p(n)}\bar{E}^p_k(f(n)\hs,n\hq)\pt^p_{\h{s}_k}v(t+\d_t,\hs,\hq+\d_q\be_k)\r\} = 0.
\end{aligned}
\]
If one assumes that there exist functions $\{\h{\mu}_k(\hs,\hq)\}_{k=1}^K$, $\{\h{\sigma}_k(\hs,\hq)\}_{k=1}^K$, such that for all $\hs\in\R^K, \hq\in[0,1]^K$,
\begin{equation}\label{eq: defofhat}
    \begin{aligned}
    &\lim_{n\to\infty}\frac{n}{f(n)}\bar{\mu}_k(f(n)\hs,n\hq)= \h{\mu}_k(\hs,\hq);\\
    &\lim_{n\to\infty}\frac{n}{f^2(n)} \bar{\sigma}^2_k(f(n)\hs,n\hq)= \h{\sigma}_k^2(\hs,\hq);\\
    &\lim_{n\to\infty}\frac{n}{f^p(n)}\bar{E}^p_k(f(n)\hs,n\hq) = \h{E}_k^p \equiv0, \quad \text{for}\quad \forall p\geq3,
\end{aligned}
\end{equation}
then as the horizon $n\to\infty$, i.e., $\d_t,\d_q\to0$, the rescaled expected cumulative reward $v(t,\hs,\hq) = \frac{1}{f(n)}w^{nt+1}(f(n)\hs,n\hq)$ satisfies the PDE
\[
\pt_tv + \max_k\{\h{\mu}_k + \pt_{\h{q}_k}v + \h{\mu}_k\pt_{\h{s}_k}v + \frac12\h{\sigma}_k^2\pt_{\h{s}_k}^2v \} = 0,
\]
which can be equivalently written as the following HJB equation:
\begin{equation}\label{eq:hjb}
    \pt_tv + \max_{\h{\bpi}(t,\hs,\hq)\in\Delta^K}\sum_{k=1}^K\l(\h{\mu}_k(t,\hs,\hq) + \pt_{\h{q}_k}v + \h{\mu}_k(\hs,\hq)\pt_{\h{s}_k}v + \frac12\h{\sigma}_k^2(\hs,\hq)\pt_{\h{s}_k}^2v \r)\h{\pi}_k = 0.
\end{equation}
Here we introduce $\h\bpi(t,\hs,\hq)$ as the feedback control, which corresponds to the policy in the bandit problem. Since the policy is a mapping from the history $H^i$ to probability measures on the action space $\Delta(\mA)$, the policy at round $i$ is described by a $K$-dimensional vector-valued function $\bpi^i(\bs^i,\bq^i)$ that satisfies $\sum_k\pi^i_k(\bs^i,\bq^i) = 1$. In the limit, the policy $\hpi(t,\hs,\hq) = \lim_{n\to\infty}\bpi^{nt+1}(f(n)\hs,n\hq)$ is a mapping from the rescaled history $(\hs,\hq)$ to the simplex $\Delta^K$ that satisfies $\sum_k\h\pi_k(t,\hs,\hq) = 1$ for $\forall (t,\hs,\hq)$.

\color{black}
We remark briefly that the selection (or even the existence) of $f(n)$ may not be obvious in the general setting described above.
In general, $f(n)$ should be thought of as describing the ``order'' or asymptotic size of the unscaled rewards in the original bandit problem with $n$ rounds. For instance, if the rewards are of constant size and observed with at least constant probability, then we will expect the cumulative reward of the original bandit problem to be linear in the time horizon, and we will have $f(n) = \Theta(n)$. Rather than describing necessary and sufficient technical conditions relating $f(n)$ to the arm reward distributions $P^{\boldsymbol{\nu}}_k$ and the posterior $\rho(\boldsymbol{\nu}|\mathbf{s}, \mathbf{q})$ (which may be intractable given the generality of the framework), in the remainder of the paper, we will show that in a wide range of concrete examples, our framework provides useful insights into the problem.
See Remark~\ref{rmk: f(n)} for further discussion of the scaling factor $f(n)$.

\color{black}
Note that the solution to \eqref{eq:hjb} is not necessarily differentiable, so we are searching for a viscosity solution instead of a classical solution \citep{evans2010partial}. One has the following guarantee on the well-posedness of the solution. 
\begin{proposition}
If  $\{\h{\mu}_k(\bs,\bq)\}_k, \{\h{\sigma}_k(\bs,\bq)\}_k$ are bounded and Lipschitz continuous in $(\bs,\bq)$, then the value function defined in \eqref{eq: def v} is the unique viscosity solution to the HJB equation \eqref{eq:hjb}.
\end{proposition}
See, e.g., \cite{nisio2015stochastic} for the proof.
\color{black}
Finding necessary and sufficient conditions on the problem primitives---specifically, the arm reward distributions $P^{\boldsymbol{\nu}}_k$, the prior $\rho$, and the scaling factor $f(n)$---under which the boundedness and Lipschitz assumptions in Proposition~1 hold is an interesting question for future work. 
In this paper, our goal is to demonstrate the useful insights which can be derived from our framework for specific bandit problems, once it has been determined that they meet these conditions. 

\paragraph{Summary} 
If the rescaled moments of all the arms satisfy \eqref{eq: defofhat} for some scaling factor $f(n)$, then the rescaled optimal expected cumulative reward $\frac{1}{f(n)}w^{nt+1}(f(n)\hs,n\hq)$ will converge to $v(t,\hs,\hq)$ as $n\to\infty$, where $v$ satisfies the HJB equation \eqref{eq:hjb} with boundary condition $v(1,\hs, \hq) = 0$.
In addition, the rescaled optimal policy $\bpi^{*,nt+1}(f(n)\hs,n\hq)$ will converge to $\hpi^*(t,\hs,\hq)$ as $n\to\infty$, where $\hpi^*(t,\hs,\hq)$ is given by
\begin{equation}\label{optimal-policy-unregu}
    \pi_k^*(t,\hs,\hq) = \l\{\begin{aligned}
       &1,\quad k = \argmax_k \l\{\h{\mu}_k(\hs,\hq) + \pt_{\h{q}_k}v + \h{\mu}_k(\hs,\hq)\pt_{\h{s}_k}v + \frac12\h{\sigma}_k^2(\hs,\hq)\pt_{\h{s}_k}^2v \r\}\\
       &0, \quad o.w.
    \end{aligned}\r.
\end{equation}
These results form the foundation for the rest of the paper.

\color{black}
\subsection{\tcb{A formal derivation} from Bayesian bandits to the optimal control problems}
If one views $v(t,\bs,\bq)$ as the optimal cumulative reward starting from time $t$
\begin{equation}\label{eq: def v}
    v(t,\hs, \hq) =\max_{\hpi(\tau)\in\Delta^K} \E\l[\int_t^1 \hat{\bm{\mu}}(\hs(\tau),\hq(\tau))\cdot\bpi(\tau) d\tau \r],
\end{equation}
then $v(t,\bs, \bq)$ is the solution to the optimal control problem \citep{evans2010partial}
\begin{equation}\label{eq:control_general}
    \begin{aligned}
    \max_{\hpi(\tau)\in\Delta^K} \ &\E\l[\int_t^1 \hat{\bm{\mu}}(\hs(\tau),\hq(\tau))\cdot\bpi(\tau) d\tau\r]\\
      s.t.\quad  &d\hlq_k(\tau) = \h{\pi}_k(\tau)d\tau,\quad 1\leq k\leq K;\\
        &d\hls_k(\tau) = \hat{\mu}_k(\hs(\tau),\hq(\tau))\h{\pi}_k(\tau)d\tau + \h{\sigma}_k(\hs(\tau),\hq(\tau))\sqrt{\h{\pi}_k(\tau)}dB_\tau,\quad 1\leq k\leq K;\\
        &\hs(t) = \hs,\quad \hq(t) = \hq.
    \end{aligned}
\end{equation}
Therefore, as the horizon $n\to\infty$, the Bayesian bandit problem also converges to the above continuous control problem. 
\tcb{In fact, one can derive the optimal control formulation above directly from the definition of the Bayesian bandit problem.}

Given a policy $\{\bpi^i\}_i$, at each round $i$, the environment $\bnu$ is sampled with probability $\rho(\bnu|\bs^i,\bq^i)$, the $k$-th arm is pulled with probability $\pi^i_k$, and the reward of the $k$-th arm follows distribution $P^\bnu_k$. 
Assume at the beginning of round $I$, $(\bs^I,\bq^I) = (\bs,\bq)$.
The goal of Bayesian bandits is to find the optimal policy $\{\bpi^{*,i}\}_i$ that maximizes the expected cumulative reward. 
More precisely, our goal is to find $\{\bpi^i\}_i$ that solves the following optimization problem:
\begin{equation}
\label{eq: dis_bb}
    \begin{aligned}
        \max_{\{\bpi^i\}_i\in\Delta^K}\quad &c_n(\{\bpi^i\}) = \E\l[\sum_{i=I}^n \int_\mE \mu_{A^i}(\bnu) \rho(\bnu|\bs^i,\bq^i)d\bnu \r]\\
        \text{where}\quad &A^i = k\quad\text{w.p. } \pi^i_k(\bs^i,\bq^i),\quad \\
        &s_k^{i+1} - s_k^{i} = \l\{ \begin{aligned}
            &0,\quad \text{if }A^i \neq k\\
            &X^i\sim P^\bnu_k, \bnu\sim\rho(\bnu|\bs^i,\bq^i),\quad \text{if  }A^i = k
        \end{aligned}
        \r.,\quad \\
        &q_k^{i+1} - q_k^i = \l\{ \begin{aligned}
            &0,\quad \text{if }A^i \neq k\\
            &1,\quad \text{if  }A^i = k
        \end{aligned}
        \r.,\\
        &(\bs^I,\bq^I) = (\bs,\bq),
    \end{aligned}
\end{equation}
where $\mu_a(\bnu)$ is the expected reward of arm $a$ defined in \eqref{eq:defmua}.

Using the same rescaling as in \eqref{eq: defofhat} and viewing the history $(\bs^i,\bq^i)$ at the discrete rounds as a function $(\hs(t),\hq(t))$ over the continuous time $t$, the differences of the rescaled cumulative reward $\hs(t)$ and the rescaled number of pulls $\hq(t)$ after one round become
\begin{equation}\label{eq:hat_diff}
    \begin{aligned}
    &\E[\hls_k(t+\dt) - \hls_k(t)] = \dt\, \frac{n}{f(n)}\pi^i_k(f(n)\hs,n\hq)\bar{\mu}_k(f(n)\hs,n\hq),\\
    &\V[\hls_k(t+\dt) - \hls_k(t)] = \dt\, \frac{n}{f^2(n)}\pi^i_k(f(n)\hs,n\hq)\bar{\sigma}^2_k(f(n)\hs,n\hq)- \l(\E[\hls_k(t+\dt) - \hls_k(t)]\r)^2,\\
    &\E[\hlq_k(t+\dt) - \hlq_k(t)] = \dt\,\pi^i_k(f(n)\hs,n\hq),\\
    &\V[\hlq_k(t+\dt) - \hlq_k(t)] = (\dt)^2\pi^i_k(f(n)\hs,n\hq) - (\dt)^2(\pi^i_k(f(n)\hs,n\hq))^2,\\
\end{aligned}
\end{equation}
where $\dt = \frac1n$, and $\bar{\mu}_k, \bar{\sigma}_k$ are defined in \eqref{eq:defofbar}. Accordingly, the rescaled objective function becomes
\begin{equation}\label{eq:hatc}
    \begin{aligned}
    \hc_n(\{\bpi^i\}) &= \frac{1}{f(n)}c_n(\{\bpi^i\})\\
    &= \dt\, \E\l[\sum_{i=I}^n\sum_{k=1}^K  \frac{n}{f(n)}\pi^i_k\l(f(n)\hs\l(\frac{i}{n}\r),n\hq\l(\frac{i}{n}\r)\r)\bar{\mu}_k\l(f(n)\hs\l(\frac{i}{n}\r),n\hq\l(\frac{i}{n}\r)\r)\r].
\end{aligned}
\end{equation}
\tcb{Formally, based on \eqref{eq:hat_diff} and \eqref{eq:hatc}, under the assumptions \eqref{eq: defofhat} on the moments, the following equations hold.
\begin{equation*}
\begin{aligned}
    &\lim_{n\to\infty}\E[\hls_k(t+\dt) - \hls_k(t)]/\dt =\h{\pi}_k(t,\hs,\hq)\hat{\mu}_k(\hs,\hq),\\
    &\lim_{n\to\infty}\V[\hls_k(t+\dt) - \hls_k(t)] /\dt= \h{\pi}_k(t,\hs,\hq)\h{\sigma}^2_k(\hs,\hq),\\
    &\lim_{n\to\infty}\E[\hlq_k(t+\dt) - \hlq_k(t)]/\dt = \h{\pi}_k(t,\hs,\hq),\\
    &\lim_{n\to\infty}\V[\hlq_k(t+\dt) - \hlq_k(t)]/\dt = 0,\\
    &\lim_{n\to\infty}\hc_n(\{\bpi^i\}) = \E\l[\int_t^1\sum_{k=1}^K\h{\pi}_k(t,\hs(\tau),\hq(\tau))\h{\mu}_k(\hs(\tau),\hq(\tau))d\tau\r].
\end{aligned}
\end{equation*}
This implies that the rescaled version of the Bayesian bandits \eqref{eq: dis_bb} will converge to the continuous optimal control problem \eqref{eq:control_general}. }

\begin{remark} \label{rmk: f(n)}
    Assumption \eqref{eq: defofhat} is critical to determine what scaling factor $f(n)$ one should use to derive a meaningful limiting HJB equation. Take the Bernoulli reward introduced in Section \ref{sec:example} for an example. Since 
    \[
        \bar{\mu}(s,q) = \bar{\s}^2(s,q) = \bar{E}^p(s,q) = \frac{\a+s}{\a+\b+q},
    \]
    one has
    \[
        \h{\mu}(\hls,\hlq) = \lim_{n\to\infty}\frac{n}{f(n)}\frac{\frac{\a}{f(n)} +\hls}{\frac{\a+\b}{n}+\hlq}, \quad \h{\s}^2(\hls,\hlq) = \lim_{n\to\infty}\frac{n}{f(n)^2}\frac{\frac{\a}{f(n)} +\hls}{\frac{\a+\b}{n}+\hlq}, 
        \]
        \[\h{E}^p(\hls,\hlq) = \lim_{n\to\infty}\frac{n}{f(n)^p}\frac{\frac{\a}{f(n)} +\hls}{\frac{\a+\b}{n}+\hlq}.
    \]
    If the initial hyperparameters $(\a,\b)$ are set such that $\lim_{n\to\infty}(\frac{\a}{f(n)},\frac{\a+\b}{n}) = (\h{\a},\h{\b})$, then the only scaling factor that will induce a meaningful HJB equation is $f(n) = O(n)$. For all $f(n) = O(n^b)$ with $b<1$, $\h{\mu}$ diverges. For the case where $b>1$, one has $\h{\mu}\equiv 0$, and consequently, the limiting HJB equation does not give any useful information on the dynamics. In this particular case, one always ends up with an HJB equation without a diffusion term ($\s \equiv 0$), and this is consistent with the nature of the Bernoulli bandits. Note that $\h{s}(t)$ is non-decreasing by definition, but when $\s>0$, there is a chance that $\h{s}(t)$ will decrease due to the nonzero $dB_\tau$ term. Therefore, any valid scaling factor must induce a deterministic optimal control problem. 
    However, in the general case, depending on the choice of scaling factor $f(n)$, the limiting optimal control problem can be stochastic or deterministic. 
\end{remark}



\subsection{Specialized limiting equations for Structured and Unstructured Bandits}\label{special eq}

In this section, we will derive the general limiting HJB equation for both unstructured and structured bandits. Then in Section~\ref{sec:3example}, we will derive the continuous limit for some specific bandit problems which are common in the literature.

\paragraph{Unstructured Bandits}
In the unstructured bandit problem, action $a\in\mA$ is completely uninformative of all other actions $b\neq a$. That is, when $a$ is played, the learner gains no information about the reward distribution of the other arms (Chapter 4.3 of \cite{lattimore2020bandit}).  Here we mainly discuss a specific kind of unstructured bandit problem in which all the arms belong to the same parametric family of distributions with different unknown parameters. For instance, the rewards of each arm may be normally distributed with an unknown mean. One will see that, in this case, the drift and diffusion terms of each arm in the limiting HJB equation have a similar form. Provided that the initial prior measure is determined by a hyperparameter $\beta$, then we claim that under certain conditions, the Bayesian bandit problem converges to the following HJB equation as the horizon $n\to \infty$:
\begin{equation}\label{eq:unstruc_hjb}
    \begin{aligned}
    &\pt_tv + \max_{\hpi(t,\hs,\hq)\in\Delta^K} \l[ \sum_{k=1}^K \l(\pt_{\hlq_k}v + \h{\mu}(\hls_k,\hlq_k,\hb_k)\pt_{\hls_k}v + \frac12\h{\sigma}(\hls_k,\hlq_k,\hb_k)^2\pt_{\hls_k}^2v \r.\r.\\
    &\quad\quad\quad\quad\quad\quad\quad\quad\quad\quad\quad\quad\quad\quad\quad\quad\quad\quad\quad\quad\quad\quad\quad\quad\quad\l.\l.+ \h{\mu}(\hls_k,\hlq_k,\hb_k)  \r) \h\pi_k(t,\hs,\hq)\r] =0,\\
    &v(1,\bs, \bq) = 0,
\end{aligned}
\end{equation}
where the functions $\h{\mu}$ and $\hat{\s}$ depend on the problem setting, and $\hb$ is the rescaled hyperparameter. The main difference between the general form \eqref{eq:control_general} and the unstructured version \eqref{eq:struc_hjb} is that the drift and diffusion functions $(\h{\mu}_k, \h{\sigma}_k)$ have the same form for each arm. Thus, the difference in drift and diffusion for each arm will be due only to differences in the arm histories $(\hls_k,\hlq_k)$ at time $t$ and the initial prior measure represented by $\h{\b}$.

\paragraph{Structured Bandits}
In a structured bandit problem, the learner can obtain information about other actions even if these actions are never played. That is, playing a particular arm $a\in\mA$ may be informative of the reward distributions of other arms $b\neq a$ (Chapter 4.3 of \cite{lattimore2020bandit}). 
Specifically, we consider the following setting. Let the action space $\{\ba_k\}_{k=1}^K = \mA\subset\R^d$ be a set of real vectors, and assume that the reward $x\in\Omega$ on the arm $\ba_k$ has density $p(x|\ba_k,\bnu)$ with an unknown parameter $\bnu\in\mE\in\R^d$. 
Then we claim that under certain conditions, the Bayesian bandit problem described above converges to the following HJB equation as the horizon $n\to\infty$:
\begin{equation}\label{eq:struc_hjb}
    \begin{aligned}
    &\pt_tv + \max_{\hpi(t,\hs,\hq)\in\Delta^K} \l[ \sum_{k=1}^K \l(\pt_{\hlq_k}v + \h{\mu}(\hs,\hq,\ba_k)\pt_{\hls_k}v + \frac12\h{\sigma}(\hs,\hq,\ba_k)^2\pt_{\hls_k}^2v \r.\r.\\
    &\quad\quad\quad\quad\quad\quad\quad\quad\quad\quad\quad\quad\quad\quad\quad\quad\quad\quad\quad\quad\quad\quad\quad\quad\quad\l.\l.+ \h{\mu}(\hs,\hq,\ba_k)  \r) \h\pi_k(t,\bs,\bq)\r] =0,\\
    &v(1,\bs, \bq) = 0,
\end{aligned}
\end{equation}
Similar to the unstructured case, the main difference between \eqref{eq:hjb} and \eqref{eq:struc_hjb} is that the drift and diffusion terms $(\h{\mu}_k,\h{\sigma}_k)$ have the same form. In the structured case, the two terms depend on the history of all the arms $(\hs,\hq)$ and the position of the arm $\ba_k$, and the difference for each arm is only due to the position  $\ba_k$.

\subsubsection{Three common bandits} \label{sec:3example}
Here we derive the HJB equation for three common examples: unstructured bandits with Bernoulli and normal arm rewards, and a structured linear bandit with normal rewards.
\paragraph{Unstructured Bernoulli rewards} 
Consider the environment class $\bnu\in[0,1]^{K}$ of Bernoulli distributions $\mB_\g$ with horizon $n$. For an environment $\bnu\in\mE = [0,1]^{K}$, the $k$-th arm follows a Bernoulli distribution taking values $\gamma(n)$ and $-\gamma(n)$ with probability $\nu_k$ and $1-\nu_k$, respectively, where $\nu_k$ is the $k$-th component of $\bnu$.  We set the initial prior measure for the $k$-th arm to be $\nu_k^1\sim\text{Beta}(\a_k(n),\b_k(n))$ for $1\leq k\leq K$.

\begin{lemma}\label{lemma: bernoulli}
Let $f(n)$ be such that there exist real numbers $(\h{\bm{\a}},\h{\bm{\b}},\h{\sigma})$ with
\[
\lim_{n\to\infty}\frac{\g(n)(\a_k(n) - \b_k(n)) }{f(n)}= \ha_k, \quad \lim_{n\to\infty}\frac{\a_k(n)+\b_k(n)}{n} = \hb_k, \quad \lim_{n\to\infty}\frac{\sqrt{n}}{f(n)}\g(n) = \h{\sigma}.
\]
Then as $n\rightarrow\infty$, the rescaled Bayesian bandit problem with scaling factor $f(n)$ converges to the HJB equation \eqref{eq:unstruc_hjb} with
\[
\h{\mu}(s,q,\ha_k,\hb_k) = \frac{\ha_k + s}{\hb_k + q}, \quad \h{\sigma}(s,q,\ha_k,\hb_k) \equiv \h{\sigma}.
\]
\end{lemma}
See Appendix~\ref{appendix: bernoulli} for the proof of the above lemma.


\paragraph{Unstructured normal rewards}
Consider the environment class $\bnu\in\R^K$ of normally distributed arm rewards $\mN_\sigma$ with horizon $n$. For environment $\bnu \in \mE = \R^K$, the rewards of the $k$-th arm follow the normal distribution $\mN(\nu_k,\sigma^2(n))$, where $\nu_k$ is the $k$-th component of $\bnu$.  
We set the initial prior measure to be $\nu_k^1\sim\mN(\a_k(n),\b_k(n)^2)$ for $1\leq k\leq K$. 
\begin{lemma}\label{lemma: normal}
Let $f(n)$ be such that there exist real numbers $(\h{\bm{\a}},\h{\bm{\b}}, \hat{\sigma})$ with
\[
\lim_{n\to\infty} \frac{\s^2(n)\a_k(n)}{f(n)\b_k^2(n)} = \ha_k, \quad \lim_{n\to\infty}\frac{\sigma^2(n)}{\b_k^2(n)n} = \hb_k,\quad \lim_{n\to\infty} \frac{\sqrt{n}}{f(n)}\sigma(n) = \h{\sigma}.
\]
Then as $n\to\infty$, the rescaled Bayesian bandit problem  with the scaling factor $f(n)$ converges to the continuous control problem \eqref{eq:struc_hjb} with
\[
\h{\mu}(s,q,\ha_k,\hb_k) = \frac{\h{\a}_k + s}{\h{\b}_k + q}, \quad \h{\sigma}(s,q,\ha_k,\hb_k) \equiv \h{\sigma}.
\]
\end{lemma}
See Appendix~\ref{appendix: normal} for the proof of the above lemma.

\paragraph{Linear bandits with normal rewards} 
Consider the case of stochastic linear bandits with horizon $n$, where the environment is encoded by a vector $\bnu \in \R^d$.  
For environment $\bnu\in \mE=\R^d$, the reward $X^i$ at round $i$ depends linearly on the chosen action $A^i\in\mA\subset\R^d$ in the following sense:
\begin{equation}\label{eq:linear}
    X^i = \la A^i, \bnu\ra + \eta^i,
\end{equation}
where $(\eta^i)_{i=1}^n$ is a sequence of independent and identically distributed normal random variables $\mN(0,\sigma^2(n))$ with given $\sigma$. 
We set the prior measure of $\bnu$ to be the normal distribution $\mN(\bba(n),\S(n))$. 

\begin{lemma}\label{lemma: linear_normal}
Let $f(n)$ be such that there exist real constants $(\h{\bba}, \h{\S})$ with
\[
\lim_{n\to\infty}\frac{\sigma^2(n)}n\S^{-1}(n) = \h{\S}^{-1},  \quad \lim_{n\to\infty}  \frac{\sigma^2(n)}{f(n)}\S^{-1}(n)\bba(n) = \h{\bba},\quad \lim_{n\to\infty} \frac{\sqrt{n}}{f(n)}\sigma(n) = \h{\sigma}.
\]
Then as $n\to\infty$, the rescaled Bayesian bandit problem with scaling factor $f(n)$ converges to the continuous control problem \eqref{eq:struc_hjb} with
\[
\begin{aligned}
    \h{\mu}(\bs, \bq,{\bf b})
    =& {\bf b}^\top\l(\h{\S}^{-1} + \sum_k     q_k\ba_k(\ba_k)^\top\r)^{-1}\l(\h{\bba} + \sum_{k=1}^Ks_k\ba_k\r), \quad 
    \h{\sigma}(\bs,\bq)  \equiv \h{\sigma}.
\end{aligned}
\]
\end{lemma}
See Appendix~\ref{appendix: linear_normal} for the proof of the above lemma. 


\section{Approximate Bayes-optimal policy}
What can one do with the limiting HJB equations or optimal control problems? In Section~\ref{sec: regu}, we propose a regularized version of the Bayesian bandit by adding a regularizer term to the objective function of the optimal control limit. In this way, the resulting optimal policy will be a stochastic policy instead of a deterministic one, which can encourage exploration and make the solution less sensitive to perturbations. In Section \ref{sec: algo}, we propose an approximate Bayes-optimal policy algorithm, which is based on the solution to the (regularized or unregularized) HJB equation. 

\subsection{Regularized Bayesian bandits}\label{sec: regu}
The optimal policy from the optimal control problem \eqref{eq:control_general} is deterministic and can be sensitive to small perturbations to the problem. To encourage robustness, one can add regularization to the objective function in \eqref{eq:control_general}, resulting in a stochastic optimal policy and encouraging more exploration and stability.  For example, entropy regularization for the policy $-\sum_k\pi_k\log\pi_k$ can be added to \eqref{eq:control_general}:
\begin{equation}\label{eq:control_regularized}
    \begin{aligned}
    \max_{\hpi(\tau)\in\Delta^K} \ &\E\l[\int_t^1 (\h{\bm{\mu}}(\hs(\tau),\hq(\tau)) - \lambda\log\hpi(\tau))\cdot\hpi(\tau)d\tau\r]\\
      s.t.\quad  &d\hlq_k(\tau) = \h\pi_k(\tau)dt,\quad 1\leq k\leq K;\\
        &d\hls_k(\tau) = \hat{\mu}_k(\hs(\tau),\hq(\tau))\h\pi_k(\tau)d\tau + \h{\sigma}_k(\hs(\tau),\hq(\tau))\sqrt{\h\pi_k(\tau)}dB_\tau,\quad 1\leq k\leq K;\\
        &\hs(t) = \hs,\quad \hq(t) = \hq.
    \end{aligned}
\end{equation}
The resulting HJB equation is \citep{evans2010partial}
\begin{equation}\label{eq: hjb1_regularized}
    \begin{aligned}
  &\pt_tv + \max_{\hpi\in\Delta^K}\l[\sum_{k=1}^{K}\l( \h{\mu}_k(\hs, \hq)\pt_{\hls_k}v+ \pt_{\hlq_k}v + \frac12\h{\sigma}_k^2(\hs, \hq)\pt_{\hls_k}^2v + \h{\mu}_k(\hs, \hq) \r.\r.\\
  &\quad\quad\quad\quad\quad\quad\quad\quad\quad\quad\quad\quad\quad\quad\quad\quad\quad\quad\quad\quad\quad\quad\quad\l.\l.- \lambda\log(\h\pi_k(t,\hs, \hq)) \r)\h\pi_k(t,\hs, \hq)\r] = 0.
\end{aligned}
\end{equation}
The maximum in the above equation can be computed explicitly \citep{ying2020note}. Let 
\[
H_k(p,m,h) = \h{\mu}_k(\bs, \bq)p + m + \frac12\h{\sigma}_k^2(\bs, \bq)h+\h{\mu}_k(\bs, \bq),
\]
then the optimal policy is
\begin{equation}\label{optimal-policy-regu}
    \h{\pi}^*_k(t,\hs,\hq) = \frac{1}{Z}\exp\l(\frac1\lambda H_k(\pt_{\hls_k}v, \pt_{\hlq_k}v, \pt_{\hls_k}^2v)\r), 
\end{equation}
with normalizing constant $Z = \sum_k\exp\l(\frac1\lambda H_k(\pt_{s_k}v, \pt_{q_k}v, \pt_{s_k}^2v)\r)$. Hence, \eqref{eq: hjb1_regularized} can be equivalently written as
\begin{equation}\label{eq:hjb-regu}
    \pt_tv + \lambda \log\l( \sum_k \exp\l( \frac1\lambda H_k(\pt_{\hls_k}v, \pt_{\hlq_k}v, \pt_{\hls_k}^2v) \r) \r) = 0.
\end{equation}
There are several potential advantages of the regularized version. First, the resulting optimal policy is always stochastic for $\lambda>0$, which will be less sensitive to perturbations compared with deterministic optimal policy. Second, regularization encourages more exploration, which helps the performance when the initial prior is significantly different from the underlying truth. Third, regularization will usually lead to a smoother solution with a differentiable policy and value function, making it easier to numerically approximate the solution.
\color{black}
For a more comprehensive introduction to the use of regularization in bandit and reinforcement learning problems, see, e.g., the tutorial by \cite{geist2020regularization}.
\color{black}

\subsection{Approximating the Bayes-optimal policy}\label{sec: algo}

Based on the limiting equation, if one can obtain the optimal policy for the HJB equation, then one can approximate the optimal Bayesian bandit policy by rescaling $(t,\hs,\hq)$ to $(i,\bs,\bq)$. This is summarized by the pseudocode in Algorithm~\ref{algo1}. 

\begin{algorithm}[h!]
\caption{Approximate Bayes-optimal policy} \label{algo1}
\begin{algorithmic}
\For{$i=1,\ldots, K$}
\State $A_i \gets i$, $s_i \gets X_i$, $q_i \gets 1$
\EndFor
\For{$i=K+1,\ldots, n$}
\State Pull $A_i \sim \bpi^i(\bs,\bq) = \h{\bpi}_k^*(\frac{i-1}{n},\frac{\bs}{f(n)}, \frac{\bq}{n})$ \Comment{$\h{\bpi}^*(t,\hs,\hq)$ given in \eqref{optimal-policy-unregu} (unregularized) or \eqref{optimal-policy-regu} (regularized)}
\State Get reward $X_i$
\State $s_{A_i}\gets s_{A_i}+X_i$, $q_{A_i}\gets q_{A_i}+1$
\EndFor
\end{algorithmic}
\end{algorithm}

\section{Solving the limiting HJB equation}
One of the difficulties of the Bayesian bandit problem is its large computational cost. The computational complexity for solving a $K$-armed bandit with horizon $n$ via backward induction is $O(n^{2K})$, which is intractable when $n$ or $K$ is large. If one can obtain the exact optimal policy for the limiting HJB equation, then one can use it to approximate the Bayes-optimal policy for the finite horizon problem with almost no additional computational cost. Section \ref{sec:exactsolu} shows one of the cases where the exact solution can be obtained. Even if the exact solution cannot be obtained directly, Section \ref{sec:numeric-scheme} shows a numerical scheme to approximate the solution. The computational cost of numerically solving the HJB equation is $O(N^{2K})$, where $N$ depends on the mesh of the scheme. This can be much more efficient than the discrete Bayesian bandit algorithm when $n$ is large and $K$ is small.

\subsection{Exact solution}\label{sec:exactsolu}
Although solving the HJB equation can also be challenging in general, it turns out that if $\mu_k(s,q) = \frac{s_k+\ha}{q_k+\hb}$, one can obtain the exact solution.  By Lemmas~\ref{lemma: bernoulli} and \ref{lemma: normal}, two common bandit problems are exactly in this form. 

\begin{theorem}\label{thm:exact}
    If the drift term in the HJB equation \eqref{eq:hjb} is
    \[\h{\mu}_k(\hs,\hq) = \frac{s_k+\ha_k}{q_k+\hb_k},\] 
    then for any constants $(\ha_k,\hb_k)$, the optimal policy for the unregularized HJB equation \eqref{eq:hjb} is 
    \begin{equation}\label{eq:pistar1}
        \h{\pi}^*_k(t,\hs,\hq) = \l\{\begin{aligned}
        &1,\quad k = \argmax_k \frac{\h{s}_k+\ha_k}{\h{q}_k+\hb_k}\\
        &0,\quad o.w.
    \end{aligned}\r.;
    \end{equation}
    and the optimal policy for the regularized HJB equation \eqref{eq: hjb1_regularized} is 
    \begin{equation}\label{eq:pistar2}
    \h{\pi}^*_k(t,\hs,\hq) \propto \exp\l(\frac{1}{\lambda}\frac{\h{s}_k+\ha_k}{\h{q}_k+\hb_k}\r).
    \end{equation}
\end{theorem}
The proof is given in Appendix~\ref{appendix: exact}. Based on the above theorem, the approximate optimal policy for the unregularized Bayesian bandit problem is given by
\begin{equation}\label{approx_unreg}
    \t{\pi}^{*,i}_k(\bs,\bq) = \h{\pi}^*_k\l(\frac{i-1}{n},\frac{\bs}{f(n)}, \frac{\bq}{n}\r) = \l\{\begin{aligned}
        &1,\quad k = \argmax_k \frac{s_k}{q_k+n\hb_k}+\frac{f(n)\ha_k}{q+n\hb_k}\\
        &0,\quad o.w.
    \end{aligned}\r.
\end{equation}
and the approximate optimal policy for the regularized Bayesian bandit problem is given by
\begin{equation}\label{approx_reg}
    \t{\pi}^{*,i}_k(\bs,\bq)= \h{\pi}^*_k\l(\frac{i-1}{n},\frac{\bs}{f(n)}, \frac{\bq}{n}\r) \propto \exp\l(\frac{n}{\lambda f(n)}\frac{s_k+f(n)\ha_k}{q_k+n\hb_k}\r).
\end{equation}

Furthermore, note that the approximate optimal policy \eqref{approx_unreg} for the unregularized Bayesian bandit is similar in form to UCB: the first term is an approximation to the empirical mean, and the second term measures the degree to which the arm has been explored. The approximate optimal policy \eqref{approx_reg} for the regularized Bayesian bandit has a form similar to the tempered greedy algorithm, where the term $\frac{s_k+f(n)\ha_k}{q_k+n\hb_k}$ is an approximation to the empirical mean and $\frac{n}{\lambda f(n)}$ adjusts the exploration rate. \tcb{When $\frac{n}{\lambda f(n)}$ is smaller (i.e., when the regularization constant $\lambda$ is larger), there is more exploration.}

\subsection{Numerical solution}\label{sec:numeric-scheme}
In the general case, an exact solution to the HJB equation \eqref{eq:hjb} is not available, so we present a numerical scheme to approximate the solution in this section. In certain cases, the numerical scheme yields the exact optimal policy and value function (see Lemma~\ref{lemma: numerics_Bayesian}). In general, when the horizon is large, the computational cost of numerically solving the PDE will be much less than that of classical Bayesian bandit algorithms while still yielding a good approximation to the optimal policy.

First, observe that one can directly compute the maximum in the HJB equation \eqref{eq:hjb}. Let
\[
k^*(t,\hs, \hq) = \text{argmax}_k\ \pt_{\hlq_k}v + \h{\mu}_k(\hs, \hq)\pt_{\hls_k}v+  \frac12\h{\sigma}_k^2(\hs, \hq)\pt^2_{\hls_k}v + \h{\mu}_k(\hs, \hq),
\]
then the optimal policy is 
\[
\pi^*_{k}(t,\hs, \hq) = \Bigg\{\begin{aligned}
&1, \quad k = k^*(t,\hs, \hq)\\
&0, \quad k \neq k^*(t,\hs, \hq)
\end{aligned}
\]
and the optimal value function satisfies
\begin{equation}\label{eq: hjb_separate}
    \begin{aligned}
&\pt_tv +  \pt_{\hlq_{k^*}}v + \h{\mu}_{k^*}(\hs, \hq)\pt_{\hls_{k^*}}v+ \frac12\h{\sigma}_{k^*}^2(\hs, \hq)\pt^2_{\hls_{k^*}}v + \h{\mu}_{k^*}(\hs, \hq)= 0.
\end{aligned}
\end{equation}
All of the above results hold in the deterministic case when $\hat{\sigma}_k\equiv 0$ for all $k$. Based on \eqref{eq: hjb_separate}, we present a finite difference method for solving the HJB equation \eqref{eq: hjb_separate}.

\paragraph{HJB equation with diffusion}
First, consider the case where $\h{\sigma}\not\equiv 0$. We discretize the time interval $[0,1]$ via the grid points $0 = t^0 < t^1 < \cdots < t^{N_t} = 1$, where $t^l = l\delta_t$ and $\d_t = 1/N_t$. We impose a cutoff on the cumulative reward in $\R^K$ so that it lies within $[-S,S]^K$, then further discretize this clipped interval into $-S \leq s_k^{\bi} \leq S$. Here $\bi\in\{\Z^K, -N_s\leq i_k\leq N_s\}$ is a $K$-dimensional index vector, and $\bs^\bi = \bi\delta_s$ with $\delta_s = S/N_s$. Finally, we discretize the number of pulls to $0\leq q^\bj_k \leq 1$, where $\bj\in\{\Z^K, 0\leq j_k\leq N_q\}$ is a $K$-dimensional index vector, and $\bq^\bj = \bj\delta_q$ with $\d_q =  1/N_q$. Observe that at each time $t^l$, one has that $0\leq q_k^j\leq t^l$ and $\sum_{k=1}^Kq_k^\bj = t^l$. We will typically set $N_q = N_t$ and $\d_q = \d_t$.

Let $\tv^{l}_{\bi,\bj}$ be the approximation for $v(t^l,\bs^\bi,\bq^\bj)$, where $\bi, \bj \in\Z^K$. The numerical approximation $\tv^l_{\bi,\bj}$ satisfies the following equation:
\begin{equation}\label{eq:nume0}
    \begin{aligned}
    &\frac1{\delta_t}(\tv^{l+1}_{\bi,\bj} - \tv^{l}_{\bi,\bj}) + \frac1{2\delta_s}\h{\mu}_{k^*}^{\bi,\bj}(\tv^{l+1}_{\bi+\be_{k^*},\bj+\be_{k^*}} - \tv^{l+1}_{\bi-\be_{k^*},\bj+\be_{k^*}})+  \frac1{\delta_q}(\tv^{l+1}_{\bi,\bj+\be_{k^*}} - \tv^{l+1}_{\bi,\bj}) 
    \\
    &+\frac1{2\delta_s^2}\h{\sigma}_{k^*}^{ \bi,\bj}(\tv^{l+1}_{\bi+\be_{k^*},\bj+\be_{k^*}} - 2\tv^{l+1}_{\bi+\be_{k^*},\bj+\be_{k^*}} + \tv^{l+1}_{\bi-\be_{k^*},\bj+\be_{k^*}})+\h{\mu}_{k^*}^{\bi,\bj} = 0,
\end{aligned}
\end{equation}
where $\be_k$ is the $k$-th standard basis vector, $\h{\mu}_{k}^{\bi,\bj} = \h{\mu}_k(\bs^\bi,\bq^\bj)$, and $\h{\sigma}_k^{\bi,\bj} = \h{\sigma}_k(\bs^\bi,\bq^\bj)$. Since the boundary conditions at the terminal time specify that $\tv^{N_t}_{\bi,\bj} = 0$, one can use the preceding equations to solve backward in time and compute all of the values $\tv^{l}_{\bi,\bj}$. 

In fact, if one defines $g^l_{\bi,\bj}$ as a vector with $k$-th element given by
\begin{equation}\label{eq:numerics_stoch1}
    \begin{aligned}
    (g^l_{\bi,\bj})_k = &\tv^{l+1}_{\bi,\bj} + \delta_t\l[\frac1{2\delta_s}\h{\mu}_k^{\bi,\bj}(\tv^{l+1}_{\bi+\be_k,\bj+\be_k} - \tv^{l+1}_{\bi-\be_k,\bj+\be_k})+  \frac1{\delta_q}(\tv^{l+1}_{\bi,\bj+\be_k} - \tv^{l+1}_{\bi,\bj}) \r.
    \\
    &\l.+\frac1{2\delta_s^2}\h{\sigma}_k^{\bi,\bj}(\tv^{l+1}_{\bi+\be_k,\bj+\be_k} - 2\tv^{l+1}_{\bi,\bj+\be_k} + \tv^{l+1}_{\bi-\be_k,\bj+\be_k})+\h{\mu}^{\bi,\bj}_k \r],
\end{aligned}
\end{equation}
then one gets an equivalent form for the numerical scheme \eqref{eq:nume0} with the correspondence given by
\begin{equation}\label{eq:numerics_stoch2}
    k^* = \text{argmax}_k (g^l_{\bi,\bj})_k, \quad \text{and}\quad \tv^l_{\bi,\bj} = (g^l_{\bi,\bj})_{k^*}, \quad 
    (\t\pi^l_{\bi,\bj})_k = \l\{
    \begin{aligned}
        &1, \quad k = k^*\\
        &0, \quad k \neq k^*
    \end{aligned}\r.
\end{equation}
Note that for this scheme to be numerically stable, $\d_t$ and $\d_s$ must satisfy the inequality
\[
\delta_t \leq \min(\h{\sigma}_{\bi,\bj}^k)^2\delta_s^2. 
\]
Since $\t{v}^l_{\bi,\bj}$ is only defined on grid points, the continuous approximated solution $\t{v}(t,\hs,\hq)$ is defined as 
\begin{equation}\label{eq:defofvhat}
    \t{v}(t,\hs,\hq) = \tv^{l}_{\bi,\bj} \quad \text{for}\quad
\l\{\begin{aligned}
    &l \d_t\leq  t <(l+1)\d_t\\
    &i_k \d_s\leq  \hls_k <(i_k+1)\d_s\\
    &j_k \d_q\leq  \hlq_k <(j_k+1)\d_q.
\end{aligned}\r.
\end{equation}

\paragraph{HJB equation without diffusion}
Second, we consider the case where $\h{\sigma} \equiv 0$. Since there is no longer a diffusion term, one should use the upwind scheme for the transport term.  Let 
\[
\h{\mu}^{\bi,\bj}_{k,+}={\text{max}}(\h{\mu}_k^{\bi,\bj},0)\,,\qquad \h{\mu}^{\bi,\bj}_{k,-}={\text{min}}(\h{\mu}_k^{\bi,\bj},0),
\]
\begin{equation}\label{eq:numerics_deter1}
\begin{aligned}
    (g^l_{\bi,\bj})_k = \tv^{l+1}_{\bi,\bj} + \delta_t&\l[\frac{\h{\mu}^{\bi,\bj}_{k,+}}{\delta_s}(\tv^{l+1}_{\bi+\be_k,\bj+\be_k} - \tv^{l+1}_{\bi,\bj+\be_k}) + \frac{\h{\mu}^{\bi,\bj}_{k,-}}{\delta_s}(\tv^{l+1}_{\bi,\bj+\be_k} - \tv^{l+1}_{\bi-\be_k,\bj+\be_k})\r.\\
    &\quad\l.+  \frac1{\delta_q}(\tv^{l+1}_{\bi,\bj+\be_k} - \tv^{l+1}_{\bi,\bj})+\h{\mu}_k^{\bi,\bj} \r].
\end{aligned}
\end{equation}
Then 
\begin{equation}\label{eq:numerics_deter2}
k^* = \text{argmax}_k (g^l_{\bi,\bj})_k \quad \text{and}\quad \tv^l_{\bi,\bj} = (g^l_{\bi,\bj})_k, \quad
    (\t\pi^l_{\bi,\bj})_k = \l\{
    \begin{aligned}
        &1, \quad k = k^*\\
        &0, \quad k \neq k^*.
    \end{aligned}\r.
\end{equation}
In this case, the stability conditions imply that
\[
\max(\h{\mu}_k^{\bi,\bj})\delta_t \leq \delta_s.
\]

\paragraph{Connection to the Bayesian bandit algorithm}
In certain cases, the numerical schemes \eqref{eq:numerics_stoch1}-\eqref{eq:numerics_stoch2} and \eqref{eq:numerics_deter1}-\eqref{eq:numerics_deter2} give the exact optimal value function for the finite horizon problem. 

\begin{lemma}\label{lemma: numerics_Bayesian}
For the Bernoulli bandits introduced in Section \ref{sec:example}, when the initial hyperparameters $(\a_k,\b_k) = (c_1n,c_2n)$ for some constants $(c_1,c_2)$, then the numerical scheme \eqref{eq:numerics_deter1}-\eqref{eq:numerics_deter2} for the limiting HJB equation based on the scaling factor $f(n) = n$ gives the exact optimal value function when $\d_t = \d_q = \d_s = \frac1n$.

Similarly, for the binomial bandits described in Section~\ref{sec:3example} with $\g$ being a constant independent of $n$, when the initial hyperparameters $(\a_k,\b_k) = (c_1n+c_2\sqrt{n},c_1-c_2\sqrt{n})$ for some constants $(c_1,c_2)$, then the numerical scheme \eqref{eq:numerics_stoch1}-\eqref{eq:numerics_stoch2} for the limiting HJB equation based on the scaling factor  $f(n) = \sqrt{n}$ gives the exact optimal value function when $\d_t = \d_q  = \frac1n, \d_s = \frac{\g}{\sqrt{n}}$.
\end{lemma}
See Appendix~\ref{appendix:numerics_Bayesian} for the proof of the above lemma.

\section{Numerical experiments}\label{sec: numerics}
\subsection{Convergence to the HJB equation}
\begin{figure}[h!]
    \centering
    \includegraphics[width=0.8\linewidth]{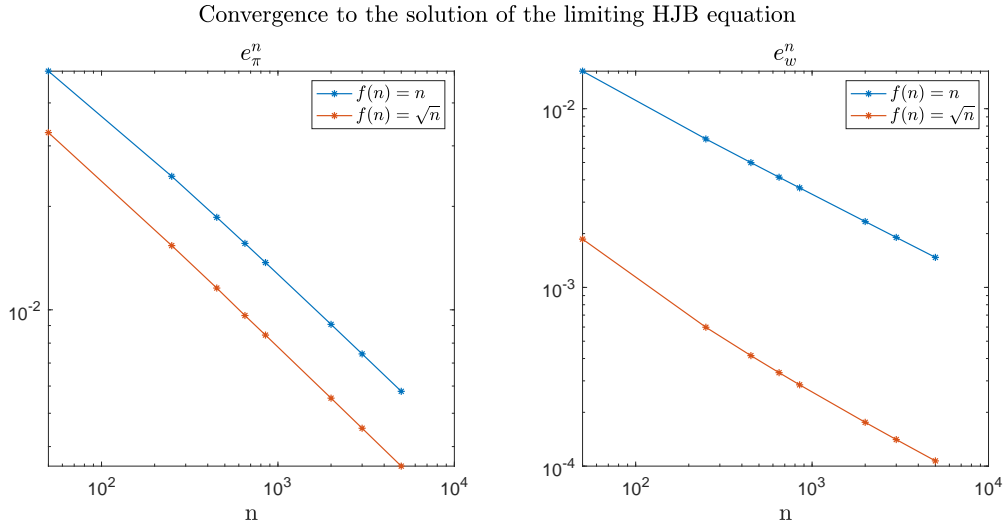}
    \caption{The above plot shows the decay of the difference between the Bayes-optimal solution and the solution to the HJB equation as $n$ increases, i.e., $e^{n}_\pi$ and $e^n_w$ defined in \eqref{error-0}. Here $f(n)$ is the scaling factor. When $f(n) = \sqrt{n}$, the resulting limit is a stochastic optimal control problem, while when $f(n) = n$, the resulting limit is a deterministic one. 
    }
    \label{fig:fig2_1}
\end{figure}

\begin{figure}[h!]
    \centering
    \includegraphics[width=0.8\linewidth]{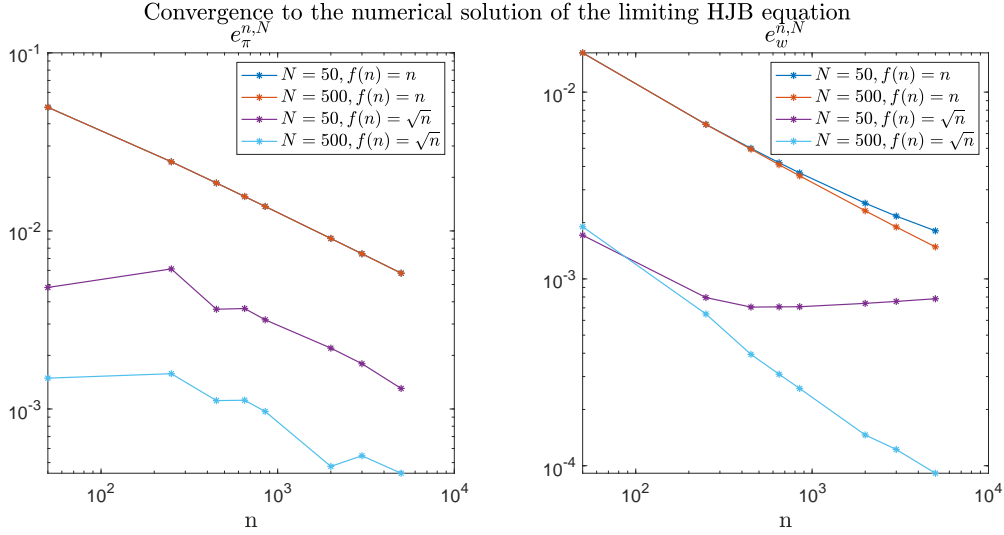}
    \caption{The above plot shows the decay of the difference between the Bayes-optimal solution and the numerical solution to the HJB equation as $n$ and $N$ increase, i.e., $e^{n,N}_\pi$ and $e^{n,N}_w$ defined in \eqref{error-1}. Here $N$ is the number of grid points when numerically solving the HJB equation, and $f(n)$ is the scaling factor. When $f(n) = \sqrt{n}$, the resulting limit is a stochastic optimal control problem, while when $f(n) = n$, the resulting limit is a deterministic one. \tcb{The curve for $N = 50,f(n) = n$ is indistinguishable from the curve for $N = 500,f(n) = n$ in $e^{n,N}_\pi$, so it is not shown on the right plot.} }
    \label{fig:fig2_2}
\end{figure}

In this section, we will show the convergence of the Bayes-optimal solution to the HJB solution as the horizon goes to infinity. 
Namely, we would like to show the differences in the optimal policy and the rescaled optimal cumulative reward
\[
\pi^{i,n}(\bs,\bq) - \h{\pi}^n\l(\frac{i-1}{n},\frac{\bs}{f(n)},\frac{\bq}{n}\r), \quad \frac{1}{f(n)}w^{i,n}(\bs,\bq) - v\l(\frac{i-1}{n},\frac{\bs}{f(n)},\frac{\bq}{n}\r)
\]
decay as the horizon $n\to\infty$, where $\pi^{i,n}(\bs,\bq), w^i(\bs,\bq)$ are obtained by backward induction and $\h{\pi}(t, \hs,\hq), v(t,\hs,\hq)$ are obtained by solving the corresponding HJB equation. We show the convergence result in Figures~\ref{fig:fig2_1} and \ref{fig:fig2_2}. Below are the details of the plots.

Consider the one-armed Bernoulli bandit problem, where the first arm has a reward $1$ with probability $\nu$ and $-1$ with probability $1-\nu$, while the second arm has a deterministic reward $\mu_2$.  In this case, given a prior measure $\nu\sim\text{Beta}(\a,\b)$, one can obtain the exact optimal policy and cumulative reward $\pi^{i,n}(s,q)$ and $w^{i,n}(s,q)$ via the equations
\[
    \pi_k^{i,n}(s,q) = \l\{\begin{aligned}
    &1,\quad \text{for}\quad k = \argmax_k\,\hat{w}^{i,n}_k(s,q)\\
    &0,\quad o.w.
    \end{aligned}\r.;
    \quad w^{i,n}(s,q) = \max_k\,\hat{w}^{i,n}_k(s,q),
\]
where
\[
\begin{aligned}
    &\hat{w}^{i,n}_1(s,q) = p(s,q)w^{i+1,n}(s+1,q+1) + (1-p(s,q)w^{i+1,n}(s-1,q+1)),\\
    &\hat{w}^{i,n}_2(s,q) = w^{i+1,n}(s,q)+\mu_2,
\end{aligned}
\]
with $w^{n+1,n}(s,q) = 0$ for all $(s,q)$ and $p(s,q) = \frac{\a + s/2 + q/2}{\a+\b +q}$.
The limiting HJB equation depends on the scaling factor $f(n)$.  By Lemma~\ref{lemma: bernoulli}, one arrives at a stochastic optimal control problem if $f(n) = \sqrt{n}$ and a deterministic one if $f(n) = n$. We will compare $\pi^{i,n},w^{i,n}$ with the limiting HJB solution for both of these scenarios.

In Figures \ref{fig:fig2_1} and \ref{fig:fig2_2}, we set $\mu_2 = 1/(3\sqrt{n})$. The hyperparameters $(\a,\b)$ for the initial prior measure are set to be $(n,n-\sqrt{n})$, which implies that in the limiting HJB equation
\[
\pt_tv + \max_{\pi(s,q)\in[0,1]}\l[\h{\mu}\pt_{\hls} v  + \pt_{\hlq} v +\frac12\h\s^2\pt_{\hls}^2v+ \h{\mu} - \mu_2\r]+\mu_2 = 0,
\]
$\h{\mu} = \frac{1+s}{2+q}, \h{\mu}_2 = 1/3, \h{\s} = 1$ for $f(n) = \sqrt{n}$, and  $\h{\mu} = \frac{s}{2+q}, \h{\mu}_2 = 0, \h{\s} = 0$ for $f(n) = n$. By Theorem~\ref{thm:exact}, one can obtain the exact optimal policy for the limiting HJB equation. In Figure~\ref{fig:fig2_1}, we plot the average difference over $i,s,q$. That is, we plot
\begin{equation}\label{error-0}
\begin{aligned}
    &e_{\pi}^{n} = \frac1Z\sum_{i,s,q}\l| \pi^{i,n}(s,q) - \h{\pi}\l(\frac{i-1}{n},\frac{s}{f(n)},\frac{q}{n}\r)\r|,\\ &e_{w}^{n} = \frac1Z\sum_{i,s,q}\l| \frac1{f(n)}w^{i,n}(s,q) - v\l(\frac{i-1}{n},\frac{s}{f(n)},\frac{q}{n}\r)\r|,    
\end{aligned}
\end{equation}
where the summation is over $i\in\{1,\cdots,n\},s\in\{-(i-1),\cdots,i-1\}, q\in\{0,\cdots,i-1\}$, and $Z$ is the number of summations. 

Next, we test the difference between the Bayes-optimal solution and the numerical solution to the HJB equation presented in Section~\ref{sec: numerics}. The following averaged differences are plotted in Figure~\ref{fig:fig2_2}:
\begin{equation}\label{error-1}
\begin{aligned}
    &e_{\pi}^{n,N} = \frac1Z\sum_{i,s,q}\l| \pi^{i,n}(s,q) - \t{\pi}^N\l(\frac{i-1}{n},\frac{s}{f(n)},\frac{q}{n}\r)\r|, \\
    &e_w^{n,N} = \frac1Z\sum_{i,s,q}\l| \frac1{f(n)}w^{i,n}(s,q) - \tv^N\l(\frac{i-1}{n},\frac{s}{f(n)},\frac{q}{n}\r)\r|,    
\end{aligned}
\end{equation}
where $\tv^N$ is obtained according to the scheme \eqref{eq:numerics_stoch1}-\eqref{eq:numerics_stoch2} with $\d_t = \d_q = N^{-1}, \d_s = N^{-1/2}$ when $f(n) = \sqrt{n}$;  and $\tv^N$ is according to the scheme \eqref{eq:numerics_deter1} - \eqref{eq:numerics_deter2} with  $\d_t = \d_q = \d_s = N^{-1}$ when $f(n) = n$. The difference $e^{n,N}_w$ based on $f(n) = \sqrt{n}$ is rescaled by $\frac{1}{\sqrt{n}}$ so that it is on the same scale as the scheme based on $f(n) = n$. Due to the different discretizations, $\tv^N\l(\frac{i-1}{n},\frac{s}{f(n)},\frac{q}{n}\r)$ is not necessarily on a grid point, so we define the continuous approximation solution $\t{v}(t,s,q)$ as in \eqref{eq:defofvhat}, which implies that 
\begin{equation}
    \t{\pi}^N\l(\frac{i-1}{n},\frac{s}{f(n)},\frac{q}{n}\r) = \t{\pi}^l_{m,j}, \quad \tv^N\l(\frac{i-1}{n},\frac{s}{f(n)},\frac{q}{n}\r) = \tv^{l}_{m, j},
\end{equation}
for 
\[
l = \l \lfloor \frac{i-1}{n\d_t}\r\rfloor, m = \l\lfloor \frac{s}{f(n)\d_s}\r\rfloor, j = \l\lfloor \frac{q}{n\d_q}\r\rfloor,
\]
where $\lfloor x \rfloor$ is the largest integer less than or equal to $x$.

Figure~\ref{fig:fig2_1} shows that the difference $e^n_\pi$ decays as $n$ increases for both scaling factors. The stochastic limit according to the scaling factor $f(n) = \sqrt{n}$ is closer to the optimal Bayesian solution compared with the deterministic limit.  Figure~\ref{fig:fig2_2} shows that the difference $e^{n,N}$ decays as $n$ and $N$ increase. Note that $e^{n,N}$ has two components: model error and numerical error. Model error decreases as the horizon $n$ increases. The numerical error decreases as the number of grid points $N$ increases.  We can see from Figure~\ref{fig:fig2_2} that when both the horizon $n$ and the number of grid points $N$ increase, the differences decrease. We observe that when $N=50$, the difference $e_v^{n,N}$ in the value functions decreases slower or does not decrease after $n$ reaches some threshold. This is because the numerical error dominates over the model error in this regime. 

\subsection{The performance of the approximate Bayes-optimal policy}\label{sec: compare} 
We compare the performance of the approximate Bayes-optimal policy (Algorithm \ref{algo1}) with Thompson sampling  and UCB  in terms of the expected regret. (See Appendix \ref{appendix: algo} for the details of Thompson sampling and UCB.)

\begin{figure}[h!]
    \centering
    \includegraphics[width=1\linewidth]{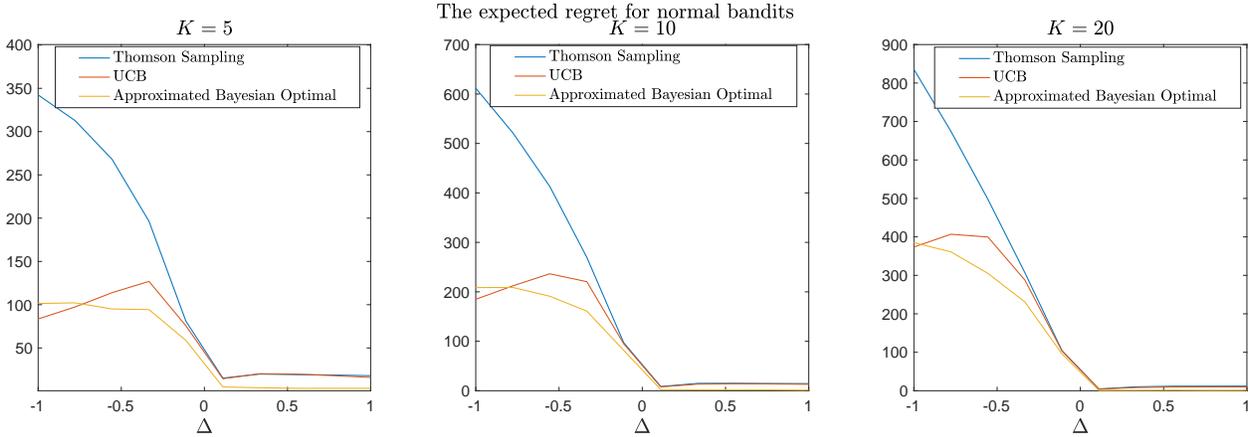}
    \caption{The above plot shows the expected regret in the $K$-armed normal bandit problem for the approximate Bayes-optimal policy, Thompson sampling, and UCB. The arm gap $\Delta\in[-1,1]$ is defined in \eqref{eg:Delta}. The left, middle, and right plots correspond to $K = 5, 10, 20$. }
    \label{fig:fig3_1}
\end{figure}

For unstructured bandits, we consider normal arm rewards, in which case the exact policy for the limiting HJB equation can be directly obtained. The performance of the three algorithms is shown in Figure~\ref{fig:fig3_1}. Below are the details of the plots.

Consider the $K$-armed normal bandit problem with $K\geq 2$. Assume that the first arm follows $r_1\sim N(\mu_1,1)$ for $\mu_1 \equiv 0$, while the $k$-th arm follows $r_k \sim N(\mu, 1)$ for $2\leq k\leq K$. Note that although $\mu_1 \equiv1$, this is unknown to us. We define
\begin{equation}\label{eg:Delta}
    \Delta = \mu - \mu_1 
\end{equation}
to be the arm gap. The horizon is set to be $n = 10^3$. For the proposed method (Algorithm \ref{algo1}), we set the scaling factor $f(n) = \sqrt{n}$, that is, the limiting optimal control problem is stochastic. The exact solution to the limiting HJB equation can be obtained by Theorem~\ref{thm:exact}. The initial prior measure for both the Bayes-optimal policy and Thompson sampling is  $\nu_k\sim N(\frac{1}{\sqrt{n}}, \frac{1}{n})$ for all $k$. This implies that the limiting HJB equation is \eqref{eq:hjb} with $\h{\mu}_k(\bs,\bq) = \frac{s_k+1}{q_k+1}$ and $\h{\s} \equiv 1$. In addition, $\delta = n^2$ for the UCB algorithm.  Figure~\ref{fig:fig3_1} shows the expected regret of the three algorithms for $\Delta\in [-1,1]$ and $K=5,10,20$. The expected regret is averaged over $10^3$ simulations. 

We can see from Figure~\ref{fig:fig3_1} that the overall performance of the approximate Bayes-optimal policy is better than the other two algorithms, especially when the prior guess is close to the underlying environment. When $\Delta$ approaches $-1$, UCB is a bit better than the approximate Bayes-optimal policy because the prior guess is significantly different from the underlying truth. However, note that as the number of arms increases, the performance is almost the same, even around $\Delta = -1$.


\begin{figure}[h!]
    \centering
    \includegraphics[width=1\linewidth]{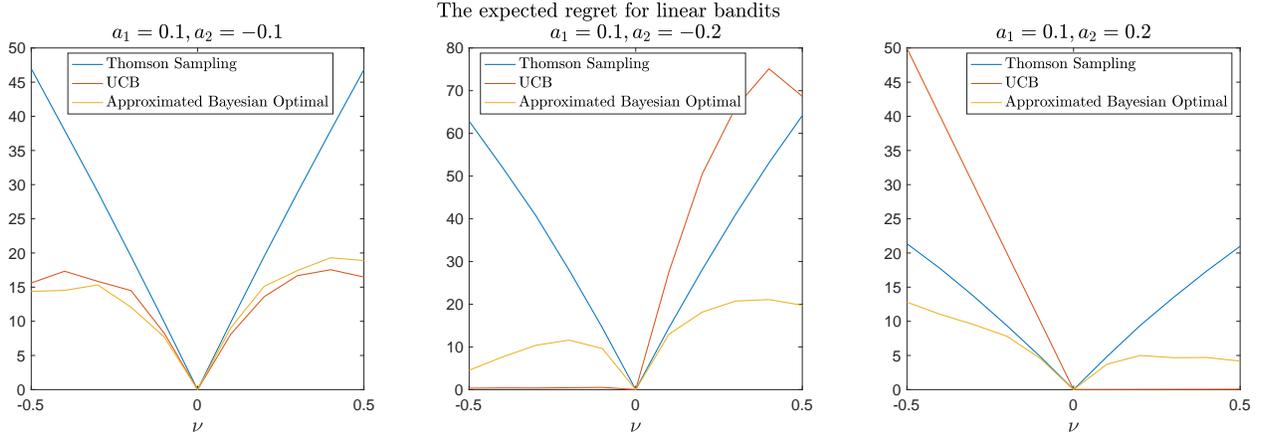}
    \caption{The above plot shows the expected regret in the $2$-armed linear bandit problem for the approximate Bayes-optimal policy, Thompson sampling, and UCB. The arm gap of the normal bandits $\Delta\in[-1,1]$ is defined in \eqref{eg:Delta}. The environment of the linear bandits $\nu \in[-1/2,1/2]$ is defined in \eqref{eq:linear}.  }
    \label{fig:fig3_3}
\end{figure}

For structured bandits, we consider the linear bandits described in Section~\ref{sec:3example}. Assume there are two arms, and the reward for arm $a_i$ follows $a_i\nu+\eta$ with unknown $\nu$ and $\eta\sim\mN(0,1)$. We set the initial measure for $\nu\sim\mN(0,\frac1n)$, and take the scaling factor $f(n) = \sqrt{n}$, then $\h\mu,\h\s$ for the limiting HJB equation can be obtained according to Lemma~\ref{lemma: linear_normal}. We solve the limiting HJB equation by the numerical scheme \eqref{eq:numerics_stoch1}-\eqref{eq:numerics_stoch2} with $\delta_t = \delta_q = \frac1N$, $\delta_s = 
\frac{1}{\sqrt{N}}$ and $N=100$.
The performance in terms of the expected regret is plotted in Figure~\ref{fig:fig3_3}, where we test three different action positions $(a_1,a_2) = (0.1,-0.1), (0.1,-0.2), (0.1,0.2)$.

Figure~\ref{fig:fig3_3} shows that the overall performance of the approximate Bayes-optimal policy is more robust than the other two methods. First, the approximate Bayes-optimal always outperforms TS. When $(a_1,a_2) = (0.1,-0.1)$, the performance of UCB and approximate Bayes-optimal policy are similar. However, for the other two cases where  $(a_1,a_2) = (0.1,-0.2), (0.1,0.2)$, UCB has much worse performance on one side, while the approximate Bayes-optimal policy has evener regret on both sides. If one measures the performance in the worst-case regret or in the averaged regret over the possible environments, the approximate Bayes-optimal policy outperforms the other two. 

\begin{figure}[h!]
    \centering
    \includegraphics[width=0.7\linewidth]{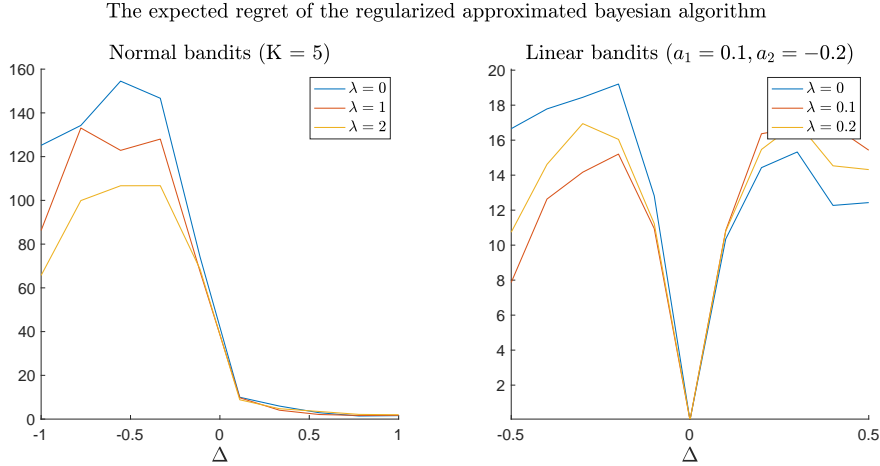}
    \caption{The above plot shows the expected regret of a $5$-armed normal bandit problem and a $2$-armed linear bandit problem for the regularized approximate Bayes-optimal policy. The environment $\nu \in[-1/2,1/2]$ is defined in \eqref{eq:linear}. The left, middle and right plots correspond to $(a_1,a_2) = (0.1,-0.1), (0.1,-0.2), (0.1,0.2)$. }
    \label{fig:fig3_2}
\end{figure}

We also show the performance of the regularized approximate Bayes-optimal policy in Figure~\ref{fig:fig3_2}. We compare the regularized Bayes-optimal policy with the unregularized version for normal bandits and linear bandits. The setting of the two bandit problem is the same as Figures \ref{fig:fig3_1} and \ref{fig:fig3_3}, but the initial Bayesian prior of the two bandits are worse (farther from the ground truth). We set $\nu_k\sim N(0.01\sqrt{n},1)$ for the normal bandits, and $\nu\sim \mN(\sqrt{n},1)$ for the linear bandits. One can see from Figure~\ref{fig:fig3_2} that the regularized version performs similarly to or better than the unregularized version when the initial prior measure is bad. 

We remark that since the solution to the regularized HJB equation is always smooth when $\lambda >0$, it is also potentially easier to break the curse of dimensionality. However, we leave the high-dimensional problem for future study. 

\section{Discussion}\label{sec:adv}
In this paper, we derived a continuous-in-time limit for the Bayesian bandit problem. We showed that the rescaled optimal cumulative reward converges to the solution of an HJB equation. We derive several benefits from the limiting PDE:
\begin{itemize}
    \item {\bf A single recipe for many Bayesian bandits.}
For most multi-armed bandit problems, the classical Bayesian bandit algorithms yield a formulation that cannot be solved accurately. Different recipes are required to approximate the value function in different settings. On the other hand, the limiting PDE gives a single, unified recipe to solve the Bayes-optimal policy.

For example, one way to solve the one-armed bandit problem with normal distributions using a Bayesian bandit algorithm is to solve the following equation backward:
\[
W^i(\mu,q) = \max\l\{W^{i+1}(\mu,q)+\mu_2, \mu + \frac{1}{2\pi}\int_\R \exp(-\frac{x}{2\sigma^2_{i}})w^{i+1}(\mu+x,q+1)dx\r\}
\]
for all $\mu\in\R$ and $i = 1, \cdots, n$ with $\sigma_i = (q^i+\sigma^{-2})^{-1}$. Since there is no closed-form solution for the integral in the above equation, one must approximate $W^i$ using piecewise quadratic functions (see e.g. Section 35.3.2 of \cite{lattimore2020bandit} for details). This results in a completely different algorithm for solving this problem compared to solving the Bernoulli reward case. However, if one instead uses the HJB equation to solve the one-armed bandit problem with normally distributed rewards, the same formulation \eqref{eq:unstruc_hjb} which applies for Bernoulli rewards also applies for normal rewards. The only modification is the different forms of $\h{\mu}$ according to Lemma \ref{lemma: normal}. 

\item {\bf Improved efficiency for large $n$.} The classical Bayesian bandit algorithm requires a computational cost of $O(n^{2K})$ to calculate the optimal policy, which can be prohibitive for large $n$.

On the other hand, as $n\to \infty$, the Bayesian bandit problem converges to the continuous HJB equation. The computational cost of solving the HJB equation is independent of the horizon, and it only depends on the numerical discretization of the PDE, which is $O(N^{2K})$. When $N\ll n$, one obtains huge computational savings by solving for the HJB value function instead. Since the accuracy of the approximation solution is $\tv - v= O(N^{-1})$, one retains an accurate approximation of the solution with much less computational cost. 

\item {\bf Improved efficiency for large $K$} For the case where the exact solution can be obtained for the limiting HJB equation as stated in Theorem \ref{thm:exact}, there is no computational cost of solving the limiting equation. In this case, even if $K$ is large, one can approximate the Bayes-optimal policy efficiently. For the case where the exact solution to the HJB equation cannot be obtained,  it may be possible to break the curse of dimensionality numerically by using a non-linear function approximation, such as a deep neural network. We leave the high-dimensional problem for future study. 
\end{itemize}

One can also extend the current framework of finite arms to infinite arms, for instance, with a continuous action space. The policy $\pi(t,\bs,\bq,\ba)$ is a probability density function such that $\int \pi(t,\bs,\bq,\ba) d\ba = 1$ for all $t,\bs, \bq$. The general limiting control problem takes the form
\begin{equation}\label{eq:lb_control_ctsaction}
    \begin{aligned}
    \max_{\int\h\pi d\ba =1}\quad &\E\l[\int_t^1\int_{\mA}\h{\mu}(\hs,\hq,\ba)\h\pi(\tau,\ba)dad\tau\r]\\
    s.t. \quad&d\hlq(\tau,\ba) = \h\pi(\tau,\ba)dt,\quad \forall \ba\in\mA;\\
    &d\hls(\tau,\ba) = \h{\mu}(\hs,\hq,\ba) \h\pi(\tau,\ba)dt + \h{\sigma}(\hs,\hq,\ba)\sqrt{\h\pi(\tau,\ba)}dB_t, \quad \forall \ba\in\mA;\\
    &\hs(t,\ba) = \bs, \quad \hq(t,\ba) = \bq.
\end{aligned}
\end{equation}
We will leave the study of the above case to future research.





\newpage
\appendix
\section*{Appendices}
\addcontentsline{toc}{section}{Appendices}
\renewcommand{\thesubsection}{\Alph{subsection}}

\subsection{Proof of Lemma \ref{lemma: bernoulli}} \label{appendix: bernoulli}
\setcounter{equation}{0}
\setcounter{theorem}{0}
\renewcommand\theequation{A.\arabic{equation}}
\renewcommand\thetheorem{A.\arabic{theorem}}
    
The expectation of the $k$-th arm in environment $\bnu$ is  $\mu_k(\bnu) = \gamma(2\nu_k - 1)$.  
At the start of the $i$-th round, the posterior distribution of $\nu_k$ is uniquely determined by the cumulative reward $s_k^i = \sum_{j=1}^{i-1}X^j\mathds{1}_{A^j = k}$ and the number of pulls $q_k^i = \sum_{j=1}^{i-1}\mathds{1}_{A^j = k}$ through round $i-1$ according to
\[
\nu_k^i \sim \text{Beta}(\a_k^i,\b_k^i), \quad \a_k^i = \a + s_k^i/(2\gamma) + q_k^i/2, \quad \b_k^i = \b - s_k^i/(2\gamma) + q_k^i/2.
\]
Hence the joint posterior $\rho^i(\bnu)$ at round $i$ is 
\begin{equation}\label{eq:bio_rho}
    \rho(\bnu|\bs^i,\bq^i) = \frac1Z\prod_{k=1}^K\nu_k^{\a_k^i-1}(1-\nu_k)^{\b^i_k-1}
\end{equation}
where $Z$ is a normalizing constant. This allows us to compute the posterior mean and variance of each arm:
\begin{equation}\label{eq:bio_bar}
    \begin{aligned}
    &\bar{\mu}_k(\bs^i,\bq^i) = \int_{[0,1]^K} \gamma(2\nu_k-1)\rho(\bnu|\bs^i,\bq^i)d\bnu = \gamma \l(2 \frac{\a_k^i}{\a_k^i + \b_k^i} - 1\r) = \gamma \frac{\a_k - \b_k+ s_k^i/\gamma}{\a_k + \b_k + q_k^i},\\
    &\bar{\sigma}^2_k(\bs^i,\bq^i) = \int_{[0,1]^K}\g^2\rho(\bnu|\bs^i,\bq^i)d\bnu = \g^2,
\end{aligned}
\end{equation}
and the higher-order moments are 
\[\bar{E}_k^p(\bs^i,\bq^i) = \int \g^p\nu_k+(-\g)^p(1-\nu_k)\rho(\bnu|\bs^i,\bq^i)d\bnu =
\l\{\begin{aligned}
    &\g^p,\quad  p\text{ is even},\\
    &\g^p\frac{\a_k - \b_k+ s_k^i/\gamma}{\a_k + \b_k + q_k^i},\quad  p\text{ is odd}.
\end{aligned}\r.\]
By the definition of $\h{\mu},\h{\s}$ in \eqref{eq: defofhat}, one has
\[
\begin{aligned}
    &\h{\mu}_k(\hs,\hq) = \lim_{n\to\infty} \frac{\frac{\gamma(\a_k - \b_k)}{f(n)} + \hls_k}{\frac{\a_k + \b_k}{n} + \hlq_k}\\
    &\h{\s}^2_k(\hs,\hq) = \lim_{n\to\infty} \frac{n}{f(n)^2}\g^2\\
    &\h{E}_k^p = \lim_{n\to\infty}\frac{n}{f(n)^p}\bar{E}_k(f(n)\bs,n\bq) = \l\{\begin{aligned}
    &\lim_{n\to\infty}\l(\frac{\g}{f(n)}\r)^{p-2}\frac{n}{f(n)^2}\g^2,\quad  p\text{ is even},\\
    &\lim_{n\to\infty}\l(\frac{\g}{f(n)}\r)^{p-1}\frac{\frac{\gamma(\a_k - \b_k)}{f(n)} + \hls_k}{\frac{\a_k + \b_k}{n} + \hlq_k},\quad  p\text{ is odd}.
\end{aligned}\r.
\end{aligned}
\]
By the assumption given in Lemma \ref{lemma: bernoulli}, one ends up with
\[
\h{\mu}(s,q,\h\a_k,\hb_k) = \lim_{n\to\infty} \frac{\h\a_k + s}{\h\b_k+q}, \quad \h{\sigma}(s,q) \equiv \h{\sigma}, \quad \h{E}^p_k\equiv0.
\]
The last equation is because $\lim_{n\to\infty}\frac{\sqrt{n}}{f(n)}\g = \h{\s}$ implies that $\lim_{n\to\infty}\frac{1}{f(n)}\g = 0$.

\paragraph{Different limiting equations due to different scaling factors}
Let us consider the case where the reward value $\g$ and $-\g$ are independent of the horizon $n$. We set $\g = 1$ and the prior hyperparameter to be $(\a,\b) = (\frac n2,\frac n2)$, i.e., the prior measure of $\nu_k\sim\text{Beta}(\frac n2,\frac n2)$. In this case, one can rescale the cumulative reward $s$ and $c_n$ by $\sqrt{n}$. By Lemma \ref{lemma: bernoulli}, one has
\[
    \h{\mu}(s,q)  = \frac{s}{1+q};\quad \hat{\sigma}(s,q) \equiv 1.
\]
One can also rescale the cumulative reward $s$ and $c_n$ by $n^{-1}$. By Lemma \ref{lemma: bernoulli}, one has
\[
    \h{\mu}(s,q)  = \frac{s}{1+q};\quad \hat{\sigma}(s,q) \equiv 0.
\]
In this case, the Bayesian bandit problem converges to a deterministic control problem in the form of \eqref{eq:unstruc_hjb} with
\begin{equation}\label{eq:bio_control2}
    \begin{aligned}
    \h{\mu}(s,q) = \frac{s}{1+q}, \quad \h{\sigma} \equiv 0.
\end{aligned}
\end{equation}
One can see from the above examples that the same Bayesian bandit problem may converge toward different control problems based on the chosen scaling.

\subsection{Proof of Lemma \ref{lemma: normal}} \label{appendix: normal}
\setcounter{equation}{0}
\setcounter{theorem}{0}
\renewcommand\theequation{B.\arabic{equation}}
\renewcommand\thetheorem{B.\arabic{theorem}}
The probability density function for the $k$-th arm in environment $\bnu$ is 
\[
P^\nu_k =\frac{1}{\sqrt{2\pi}\sigma}e^{-\frac{|x-\nu_k|^2}{2\sigma^2}}.
\]
(Note that here $\pi$ denotes the constant 3.14159... rather than the policy.)
Thus the expected reward for the $k$-th arm in environment $\bnu$ is $\mu_k(\bnu) = \nu_k$. At the start of the $i$-th round, the posterior distribution of $\nu_k$ is uniquely determined by the cumulative reward $s_k^i$ and the number of pulls $q_k^i$ through round $i-1$ according to
\[
\nu_k^i \sim \mN(\a_k^i,(\b_k^i)^2), \quad \a_k^i = \frac{\a\b^{-2} + s^i_k\sigma^{-2}}{q_k^i\sigma^{-2} + \b^{-2}}, \quad (\b_k^i)^2 = \frac{1}{q_k^i\sigma^{-2} + \b^{-2}}.
\]
Hence the joint posterior $\rho^i(\bnu)$ at round $i$ is
\[
\rho^i(\bnu) = \rho(\bnu|\bs^i,\bq^i) =\prod_{k=1}^K\frac{1}{\sqrt{2\pi}\b^i_k}e^{-\frac{|\nu_k - \a_k^i|^2}{2(\b^i_k)^2}}.
\] 
From this, it follows that
\[
\begin{aligned}
&\bar{\mu}_k(\bs^i,\bq^i) = \int_{\R^K} \nu_k\rho^i(\bnu)d\bnu = \a_k^i= \frac{\a\b^{-2} + s^i_k\sigma^{-2}}{q_k^i\sigma^{-2} + \b^{-2}},\\
&\bar{\sigma}^2_k(\bs^i,\bq^i) = \int_{\R^K} (\sigma^2+\nu_k^2)\rho^i(\bnu)d\bnu = \sigma^2 + (\b_k^i)^2 + (\a_k^i)^2=\sigma^2  + \frac{1}{q^i_k\sigma^{-2} + \b^{-2}} + (\bar{\mu}_k)^2,
\end{aligned}
\]
Inserting $\bar{\mu}_k$ and $\bar{\sigma}_k$ into \eqref{eq: defofhat} yields
\[
\begin{aligned}
    \h{\mu}_k(\hs,\hq) =& \lim_{n\to\infty} \frac{n}{f(n)} \frac{\a\b^{-2} + f(n) \hls_k\sigma^{-2}}{n\hlq_k\sigma^{-2} + \b^{-2}}= \lim_{n\to\infty}\frac{\hls_k+\frac{\sigma^2\a}{f(n)\b^2}}{\hlq_k + \frac{\sigma^2}{n\b^2}}.\\
    (\h{\sigma}_k(\hs,\hq))^2 =&\lim_{n\to\infty}\frac{n}{f(n)^2}\sigma^2 + \frac{n}{f(n)^2(n\sigma^{-2}\hlq_k+\b^{-2})} + \frac{1}{n}\l(\frac{n}{f(n)}\bar{\mu}_k\r)^2\\
=&\lim_{n\to\infty}\frac{n}{f(n)^2}\sigma^2\l(1 + \frac{1}{n}\frac{1}{\hlq_k+\frac{\sigma^2}{n\b^{2}}}\r).\\
\end{aligned}
\]
By the condition in Lemma \ref{lemma: normal}, one ends up with
\[
\h{\mu}(s,q,\h\a_k,\hb_k) =  \frac{\h\a_k + s}{\h\b_k^{-2}+q}, \quad \h{\sigma}(s,q) \equiv \h{\sigma}.
\]
For the higher-order moments, one can write the moments of the normal distribution in the following form,
\[
    \int x^p P^\nu_k(x) dx = \sum_{j=0}^{\lfloor p/2 \rfloor}C(p,j)\nu_k^{p-2j}\s^{2j}
\]
where $C(p,j)$ is a constant depends on $p,j$. Then one has
\[
\begin{aligned}
&\bar{E}^p(\bs^i,\bq^i) = \int \int x^p P^\nu_k(x) dx\rho(\nu_k)d\nu_k = \sum_{j=0}^{\lfloor p/2 \rfloor}C(p,j)\s^{2j}\int\nu_k^{p-2j}\rho(\nu_k)d\nu_k\\
=&\sum_{j=0}^{\lfloor p/2 \rfloor} \sum_{l=0}^{\lfloor (p-2j)/2 \rfloor} C(p,j)C(p-2j,l)\s^{2j}(\b^i_k)^{2l}(\a_k^i)^{p-2j-2l}\\
=&\sum_{j=0}^{\lfloor p/2 \rfloor} \sum_{l=0}^{\lfloor (p-2j)/2 \rfloor} C(p,j)C(p-2j,l)\s^{2j}(\b^i_k)^{2l}(\bar{\mu}_k)^{p-2j-2l}
\end{aligned}
\]
Therefore, after rescaling, one has
\[
\begin{aligned}
&\h{E}_k^p(\hs,\hq) =\lim_{n\to\infty} \sum_{j=0}^{\lfloor p/2 \rfloor} \sum_{l=0}^{\lfloor (p-2j)/2 \rfloor}C(p,j)C(p-2j,l)n^{1+j-p}\l(\frac{n}{f(n)^2}\s^2\r)^{j+l}\l(\frac{1}{\hlq_k+\frac{\s^2}{n\b^2}}\r)^l\l(\frac{n}{f(n)}\bar{\mu}_k\r)^{p-2j-2l}
\end{aligned}
\]
Since 
\[
n^{1+j-p}\leq n^{-1}, \quad \lim_{n\to\infty}\frac{n}{f(n)^2}\s^2 = \h\s, \quad \lim_{n\to\infty}\frac{1}{\hlq_k+\frac{\s^2}{n\b^2}} = \frac{1}{\hlq_k+\h\b^{-2}},\quad \lim_{n\to\infty}\frac{n}{f(n)}\bar{\mu}_k = \h\mu,
\]
one has,
\[
\h{E}_k^p(\hs,\hq) \equiv 0
\]

\subsection{Proof of Lemma \ref{lemma: linear_normal}} \label{appendix: linear_normal}
\setcounter{equation}{0}
\setcounter{theorem}{0}
\renewcommand\theequation{C.\arabic{equation}}
\renewcommand\thetheorem{C.\arabic{theorem}}

The expectation of the reward at round $i$ in environment $\bnu$ is $\la A^i, \bnu\ra$, and the probability density function for $k$-th arm is
\[
P^\bnu_k = \frac{1}{\sqrt{2\pi}\sigma}e^{-\frac{|x - (\ba_k)^\top\bnu|}{2\sigma^2}}.
\]
Then, the posterior distribution of $\bnu$ at round $i$ is uniquely determined by the cumulative reward $(\bs^i,\bq^i)$ up to time $i-1$,
\[
\bnu^i \sim \mN(\bba^i,\S^i), \quad \S^i = (\S^{-1} + \sigma^{-2}\sum_{k=1}^K q_k^i\ba_k (\ba_k)^\top)^{-1},\quad \bba^i = \S^i(\S^{-1}\bba + \sigma^{-2}\sum_{k=1}^Ks^i_k\ba_k).
\]
Note that the $\ba_k\in\R^d$ in the above equation represents the action value of the $k$-th arm.
Hence the posterior measure of $\bnu^i$ at round $i$ is
\[
\rho^i(\bnu) = \frac{1}{\sqrt{(2\pi)^d|\S^i|}}\exp\l(-\frac12(\bnu - \bba^i)^\top\S^i(\bnu - \bba^i)\r).
\]
Therefore, one obtains 
\begin{equation}\label{eq:lb_bar}
    \begin{aligned}
\bar{\mu}_k(\bs^i,\bq^i) =& \int_{\R^K} (\ba_k)^\top\bnu\rho^i(\bnu)d\bnu = (\ba_k)^\top\bba^i\\
=&(\ba_k)^\top(\S^{-1} + \sigma^{-2}\sum_{j=1}^Kq^i_j \ba_j(\ba_j)^\top)^{-1} (\S^{-1}\bba + \sigma^{-2}\sum_{j=1}^Ks^i_j\ba_j),\\
\bar{\sigma}^2_k(\bs^i,\bq^i) =& \int_{\R^K} (\sigma^2+(\ba_k)^\top\bnu\bnu^\top(\ba_k)^\top)\rho^i(\bnu)d\bnu = \sigma^2 +  (\ba_k)^\top(\S^i + \bba^i(\bba^i)^\top)\ba_k\\
=&\sigma^2  + (\ba_k)^\top(\S^{-1} + \sigma^{-2}\sum_{j=1}^Kq^i_j\ba_j(\ba_j)^\top)^{-1}\ba_k + (\bar{\mu}_k)^2.
\end{aligned}
\end{equation}
Inserting $\bar{\mu}_k$ and $\bar{\sigma}_k$ into \eqref{eq: defofhat} yields
\[
\begin{aligned}
    &\h{\mu}_k(\hs,\hq) = \lim_{n\to\infty}  (\ba_k)^\top\l(\frac{\sigma^2}n\S^{-1} + \sum_{j=1}^K\hlq^i_j \ba_j(\ba_j)^\top\r)^{-1} \l(\frac{\sigma^2}{f(n)}\S^{-1}\bba + \sum_{j=1}^K\hls^i_j\ba_j\r).\\
    &(\h{\sigma}_k(\hs,\hq))^2 =\lim_{n\to\infty}\l(\frac{n}{f(n)^2}\s^2\r)\l(1  + (\ba_k)^\top\l(\frac{\s^2}n\S^{-1} + \sum_{j=1}^K\hlq^i_j\ba_j(\ba_j)^\top\r)^{-1}\ba_k\r) + \frac1{n}\l(\frac{n}{f(n)}\bar{\mu}_k\r)^2.\\
\end{aligned}
\]
By the condition in Lemma \ref{lemma: normal}, one ends up with
\[
\h{\mu}(\bs,\bq,{\bf b}) = {\bf b}^\top(\h{\S}^{-1} + \sum_k     q_k\ba_k(\ba_k)^\top)^{-1}(\h{\bba} + \sum_{k=1}^Ks_k\ba_k), \quad \h{\sigma}(\bs,\bq) \equiv \h{\sigma}.
\]
For the higher-order moments, since
\[
\begin{aligned}
    &\S(\hs,\hq) = \s^2n^{-1}\l(\frac{\s^2}{n}\S^{-1}+\sum_k\hlq_k\ba_k(\ba_k)^\top\r)^{-1}\\
    &\bba^i(\hs,\hq) = \frac{f(n)}{n}\l(\S^{-1}+\sum_k\hlq_k\ba_k(\ba_k)^\top\r)^{-1}\l(\frac{\s^2\S^{-1}\bba}{f(n)}+\sum_k\hls_k\ba_k\r),
\end{aligned}
\]
Therefore, similar to the normal bandit problem, the higher-order moments are in the following order:
\[
\begin{aligned}
    \h{E}_k(\hs,\hq) = &\lim_{n\to\infty}\frac{n}{f(n)^p}\sum_{j=0}^{\lfloor p/2\rfloor}\sum_{l = 0}^{\lfloor (p-2j)/2\rfloor}O(\s^{2j})O(\S^l)O(\bba^{p-2j-2l})\\
    =&\lim_{n\to\infty}\sum_{j=0}^{\lfloor p/2\rfloor}\sum_{l = 0}^{\lfloor (p-2j)/2\rfloor} O\l(\frac{n}{f(n)^p}\r)O(\s^{2j})O(\s^{2l}n^{-l})O\l(\frac{f(n)^{p-2j-2l}}{n^{p-2j-2l}}\r)\\
    =&\lim_{n\to\infty}\sum_{j=0}^{\lfloor p/2\rfloor}\sum_{l = 0}^{\lfloor (p-2j)/2\rfloor} O\l(\l( \frac{\sqrt{n}\s}{f(n)}\r)^{2j+2l} {n^{1+j-p}}\r) = 0
\end{aligned}
\]
where the last equality is because $1_j-p\leq -1$ for all $p\geq3$ and $\lim_{n\to\infty}\frac{\sqrt{n}\s}{f(n)} = \h\s$.


\subsection{Proof of Theorem \ref{thm:exact}}\label{appendix: exact}
\setcounter{equation}{0}
\setcounter{theorem}{0}
\renewcommand\theequation{D.\arabic{equation}}
\renewcommand\thetheorem{D.\arabic{theorem}}
Look at the optimal control problem \eqref{eq:control_regularized} with $\h\s\equiv 0$, when $\h\mu_k = \frac{\hls_k+\a_k}{\hlq_k+\b_k}$, note that
\[
\frac{d}{d\tau}\h\mu_k = \frac{1}{\hlq_k+\b_k}\frac{d}{d\tau}\hls_k - \frac{\hls_k+\a_k}{(\hlq_k+\b_k)^2}\frac{d}{d\tau}\hlq_k = 0 \quad \text{for}\quad \forall \h\pi(\tau).
\]
This implies that $\h\mu_k(\tau) \equiv \h\mu_k(t) $ for $\forall k$. Therefore, the objective function becomes
\begin{equation}\label{eq:obj1}
    \int_t^1 (\h{\bm{\mu}}(t)-\lambda\bpi(\tau))\cdot\bpi(\tau) d\tau
\end{equation}
Since $\h{\bm{\mu}}(t)$ is a constant, the above objective function will be maximized at $\pi^*$ given in \eqref{eq:pistar1} and \eqref{eq:pistar2} for the unregularized version, i.e., $\lambda =0$ and the regularized version, i.e.,  $\lambda >0$.

For the stochastic case, Note that 
\begin{equation}
    d\E[\h\mu_k] = \E\left[\frac{1}{\hlq_k+\b_k}d\hls_k - \frac{\hls_k+\a_k}{(\hlq_k+\b_k)^2}d\hlq_k\right] = 0+\E[\h{\s}(\hs,\hq)\sqrt{\pi(\hs,\hq)}dB_t] =0\quad \text{for}\quad \forall \h\pi(\tau).
\end{equation}
so the objective function for the stochastic case is the same as \eqref{eq:obj1}, which results in the same optimal policy. 


\subsection{Proof of Lemma \ref{lemma: numerics_Bayesian}} \label{appendix:numerics_Bayesian}
\setcounter{equation}{0}
\setcounter{theorem}{0}
\renewcommand\theequation{E.\arabic{equation}}
\renewcommand\thetheorem{E.\arabic{theorem}}

Consider $k$-armed Bernoulli bandit, where the $k$-th arm gives reward $1$ with probability $\nu_k$ and $0$ with probability $1-\nu_k$. The initial prior measure of $\nu_k$ follows $\text{Beta}(\a_k,\b_l)$. Then the rescaled optimal cumulative reward
\[
v^i(\hs,\hq) = \frac{1}{n}w^i(\bs,\bq)
\]
with scaling factor $f(n) = n$ satisfies
\begin{equation}\label{eq:v-dis-deter}
    v^i(\hs,\hq) = \max_{k}\l\{ \frac1n\t{p}_k(\hs,\hq)+\t{p}_k(\hls,\hlq)v^{i+1}(\hs+\frac1n\be_k,\hlq+\frac1n\be_k)+(1-\t{p}_k(\hs,\hq))v^{i+1}(\hs,\hq+\frac1n\be_k) \r\}
\end{equation}
where 
\[
\t{p}_k(\hs,\hq) = \frac{n^{-1}\a_k + s_k}{n^{-1}(\a_k+\b_k)+q_k}.
\]
The corresponding HJB equation under the scaling factor $f(n) = n$ is
\[
\pt_tv + \max_{\hpi(t,\hs,\hq)\in\Delta^K}\l(\h{\mu}_k + \h{\mu}_k\pt_{s_k}v + \pt_{q_k}v\r)\pi_k  = 0, \quad \quad v(1,\hls,\hlq) = 0.
\]
where 
\[
\h{\mu}_k = \frac{\h\a_k+s_k}{\h\b_k+q_k}, \quad\text{with}\quad\h\a_k = \lim_{n\to\infty}n^{-1}\a_k, \quad\h\b_k = \lim_{n\to\infty}n^{-1}(\a_k+\b_k).
\]
Since $\h\mu \geq 0$, applying the numerical scheme \eqref{eq:numerics_deter1}-\eqref{eq:numerics_deter2} with $\d_t = \d_q = \d_s = 1/n$ gives ,
\begin{equation}\label{eq:v-num-deter}
    \begin{aligned}
    \t{v}(\frac{l}n,\hs,\hq) =\max_k&\l\{\t{v}(\frac{l+1}n,\hs,\hq) + \frac1n\l[ \frac{\h{\mu}_k(\hs,\hq)}{n^{-1}}\l(\t{v}(\frac{l+1}n,\hs+\frac1n\be_k,\hq+\frac1n\be_k)  \r.\r.\r.\\
    &\l.-\t{v}(\frac{l+1}n,\hs,\hq+\frac1n\be_k)\r) + \frac{1}{n^{-1}}\l(\t{v}(\frac{l+1}n,\hs,\hq+\frac1n\be_k) - \t{v}(\frac{l+1}n,\hs,\hq)\r) \\
    &\l.\l.+ \h{\mu}_k(\hs,\hq) \r] \r\}
\end{aligned}
\end{equation}
By comparing \eqref{eq:v-dis-deter} and \eqref{eq:v-num-deter}, one can see if and only if 
\[
\h\mu_k(\hs,\hq) = \t{p}_k(\hs,\hq),
\]
The numerical scheme is equivalent to the exact Bayes-optimal algorithm. The above condition holds if and only if 
\[
(\a_k,\b_k) = (c_1n,c_2n)
\]
for any constants $(c_1,c_2)$, which completes the proof for the first part of the Lemma \ref{lemma: numerics_Bayesian}.

Consider the binomial bandits described in Section \ref{sec:3example}. First, the optimal cumulative reward satisfies,
\[
\begin{aligned}
 w^i(\bs,\bq) = \max_k&\l\{\g(2p_k(\bs,\bq)-1) + p_k(\bs,\bq)w^{i+1}(\bs+\g\be_k,\bq+\be_k)\r.\\
 &\ \ \l.+(1-p_k(\bs,\bq))w^{i+1}(\bs-\g\be_k,\bq+\be_k)\r\}    
\end{aligned}
\]
with $p_k(\bs,\bq) = \frac{\a_k+s_k/(2\g)+q_k/2}{\a_k+\b_k+q_k}$.
Then the rescaled optimal cumulative reward
\[
v^i(\hs,\hq) = \frac{1}{\sqrt{n}}w^i(\bs,\bq)
\]
with scaling factor $f(n) = \sqrt{n}$ satisfies
\[
    \begin{aligned}
    v^i(\hs,\hq) = \max_k&\l\{n^{-1/2}\g(2\h{p}_k(\hs,\hq)-1)+\h{p}_k(\hs,\hq)v^{i+1}(\hs+\g n^{-1/2}\be_k,\hq+n^{-1}\be_k)\r.\\
    &\quad\quad\quad\quad\quad\quad\quad\quad\quad\l. +(1-\h{p}_k(\hs,\hq))v^{i+1}(\hs-\g n^{-1/2}\be_k,\hq+n^{-1}\be_k)\r\},
\end{aligned}
\]
where 
\[
    2\h{p}_k(\hs,\hq)-1 = \frac{1}{\sqrt{n}\g}\frac{\frac{\g(\a_k-\b_k)}{\sqrt{n}}+\hls_k }{\frac{\a_k+\b_k}n+\hlq_k}
\]
By letting $\t{\mu}(\hs,\hq) = \sqrt{n}\g(2\h{p}_k(\hs,\hq)-1)$, then 
\begin{equation}\label{eq:v-dis-stoch}
    \begin{aligned}
    v^i(\hs,\hq) = \max_k&\l\{n^{-1}\t{\mu}(\hs,\hq)+\frac12\l(\frac{\t{\mu}(\hs,\hq)}{\sqrt{n}\g} + 1\r)v^{i+1}(\hs+\g n^{-1/2}\be_k,\hq+n^{-1}\be_k)\r.\\
    &\l. +\frac12\l(1-\frac{\t{\mu}(\hs,\hq)}{\sqrt{n}\g} \r)v^{i+1}(\hs-\g n^{-1/2}\be_k,\hq+n^{-1}\be_k)\r\},
\end{aligned}
\end{equation}
On the other hand, the limiting HJB equation for $v^i(\hs,\hq)$ is
\[
\pt_tv + \max_{\hpi(t,\hs,\hq)\in\Delta^K}\l(\h{\mu}_k + \h{\mu}_k\pt_{s_k}v + \pt_{q_k}v + \frac12\h\s^2\pt_{s_k}^2v\r)\pi_k  = 0, \quad \quad v(1,\hls,\hlq) = 0,
\]
where \[
\h\mu(\hs,\hq) = \frac{\h\a_k+\hls_k}{\h\b_k+\hlq_k}, \quad \h{\sigma} = \g, \quad \text{with}\quad\ha_k = \lim_{n\to\infty}\frac{\g(\a_k - \b_k) }{\sqrt{n}}, \quad \hb_k = \lim_{n\to\infty}\frac{\a_k+\b_k}{n}.
\]
By comparing \eqref{eq:v-dis-deter} and \eqref{eq:v-num-deter}, one can see if and only if 
\[
\h\mu_k(\hs,\hq) = \t{\mu}_k(\hs,\hq),
\]
The numerical scheme is equivalent to the exact Bayes-optimal algorithm. The above condition holds if and only if 
\[
(\a_k,\b_k) = (c_1n+c_2\sqrt{n},c_1-c_2\sqrt{n})
\]
for any constants $(c_1,c_2)$, which completes the proof for the second part of the Lemma \ref{lemma: numerics_Bayesian}.

\subsection{Detailed algorithms}\label{appendix: algo}
\setcounter{equation}{0}
\setcounter{theorem}{0}
\renewcommand\theequation{F.\arabic{equation}}
\renewcommand\thetheorem{F.\arabic{theorem}}

\begin{algorithm}[h!]
\caption{Thompson Sampling for unstructured bandits} \label{alg:thompsonunstruc}
\begin{algorithmic}
\Require Input: $n, (\bs,\bq) = 0, \rho(\bnu|\beta)$
\For{$i=1,\ldots, K$}
\State $A_i \gets i$, $s_i \gets X_i$, $q_i \gets 1$
\EndFor
\For{$i=K+1,\ldots, n$}
\State Update $\rho(\bnu|\bs,\bq)$ according to Bayesian rule
\State Sample $\bnu^i$ according to the probability distribution $\rho(\bnu|\bs,\bq)$
\State $A_i \gets \argmax_k \mu_k(\bnu^i)$, $s_k \gets s_k+X_i$, $q_k \gets q_k+1$
\EndFor
\end{algorithmic}
\end{algorithm}

\begin{algorithm}[h!]
\caption{UCB for unstructured bandits} \label{alg:ucbunstruc}
\begin{algorithmic}
\Require Input: $n, (\bs,\bq) = 0, \delta$
\For{$i=1,\ldots, K$}
\State $A_i \gets i$, $s_i \gets X_i$, $q_i \gets 1$
\EndFor
\For{$i=K+1,\ldots, n$}
\State Update $\hat{\mu}_k \gets \frac{s_k}{q_k} + \sqrt{\frac{2\log(\delta)}{q_k}}$
\State $A_i \gets \argmax_k \hat{\mu}_k$, $s_k \gets s_k+X_i$, $q_k \gets q_k+1$
\EndFor
\end{algorithmic}
\end{algorithm}

\begin{algorithm}[h!]
\caption{Thompson Sampling for linear bandits} \label{alg:thompsonlin}
\begin{algorithmic}
\Require Input: $n, \ba_k, (\bs,\bq) = 0, \rho(\bnu)\sim\mN(\mu,\S)$
\For{$i=1,\ldots, K$}
\State $A_i \gets i$, $s_i \gets X_i$, $q_i \gets 1$
\EndFor
\For{$i=K+1,\ldots, n$}
\State Sample $\bnu$ according to the probability distribution $\rho(\bnu)$
\State $A^i = \argmax_k \ba_k^\top\bnu$
\State $x = \ba_{A^i}$, $y = \bnu*x+\eta$, where $\eta\sim\mN(0,\s^2)$
\State $\mu = (\S^{-1}+\s^{-2}xx^\top)^{-1}(\S^{-1}\mu+\s^{-2}yx)$, $\S = (\S^{-1}+\s^{-2}xx^\top)^{-1}$
\EndFor
\end{algorithmic}
\end{algorithm}

\begin{algorithm}[h!]
\caption{UCB for linear bandits} \label{alg:ucblin}
\begin{algorithmic}
\Require Input: $V=v_0, W = 0, \h\bnu = 0, \lambda = 0.1$
\For{$i=1,\ldots, K$}
\State $A_i \gets i$, $s_i \gets X_i$, $q_i \gets 1$
\EndFor
\For{$i=K+1,\ldots, n$}
\State $\b \gets \sqrt{\lambda} + \sqrt{2\log(n^2)+\log(1+(i-1)/\lambda)}$
\State $A^i \gets \argmax_k \ba_k\h{\bnu}+\b*\sqrt{\ba_k^2/V}$
\State $x \gets \ba_{A^i}$, $y \gets \bnu*x+\eta$, where $\eta\sim\mN(0,\s^2)$
\State $V \gets V+x^2, W \gets W+xy, \h{\bnu} \gets \frac{W}{V}$
\EndFor
\end{algorithmic}
\end{algorithm}

\end{document}